\documentclass{article}
\usepackage{graphicx} 
\usepackage{mathrsfs}
\usepackage{mathrsfs}
\usepackage{amsfonts}
\usepackage{amsmath}
\usepackage{amssymb}
\usepackage{float}
\usepackage{multirow}
\usepackage[all]{xy}
\usepackage{graphicx}
\usepackage{xcolor}
\usepackage{subcaption}
\usepackage[numbers]{natbib}

\allowdisplaybreaks

\newenvironment{keywords}{%
  \par\medskip\noindent
  \small
  \textbf{Keywords:}\ \ignorespaces
}{%
  \par\medskip
}
\newtheorem{exmp}{Example}
\newtheorem{theorem}{Theorem}
\newtheorem{lemma}{Lemma}
\newtheorem{corollary}{Corollary}

\title{$Q$ statistics in data depth: fundamental theory revisited and variants}
\author{Min Gao \\ 
School of Big Data and Statistics, 
Anhui University, Hefei, PR China 230601;\\ Department of Computer Science, Mathematics, Physics and Statistics, \\ University of British Columbia, Kelowna, Canada V1V 1V7
\and Yiting Chen \\ 
Department of Computer Science, Mathematics, Physics and Statistics, \\ University of British Columbia, Kelowna, Canada V1V 1V7
\and Xiaoping Shi \\
Department of Computer Science, Mathematics, Physics and Statistics, \\ University of British Columbia, Kelowna, Canada V1V 1V7
\and Wenzhi Yang \\ 
School of Big Data and Statistics, 
Anhui University, Hefei, PR China 230601
}

\begin{document}
\maketitle
\begin{abstract} 
Recently, data depth has been widely used to rank multivariate data.
The study of the depth-based $Q$ statistic,  originally proposed by Liu and Singh (1993), has become increasingly popular when it    can be used as a quality index to differentiate between two samples. Based on the existing theoretical foundations, more and more variants have been developed for increasing power in the two sample test.
However, the asymptotic expansion of the $Q$ statistic in the important foundation work of Zuo and He (2006) currently has an optimal rate $m^{-3/4}$ slower than the target $m^{-1}$, leading to limitations in higher-order expansions for developing more powerful tests.  
 We revisit the existing assumptions and add two new plausible  assumptions to obtain the target rate by applying a new proof method based on the Hoeffding decomposition and the Cox-Reid expansion.  
 The aim of this paper is to rekindle interest in asymptotic data depth theory, to place Q-statistical inference on a firmer theoretical basis, to show its variants in current research, to open the door to the development of new theories for further variants requiring higher-order  expansions, and to explore more of its potential applications.

\end{abstract}

\section{Introduction}
The two-sample test is a basic tool for comparing two groups or populations and is widely used in scientific research, engineering, and medical practice.
Common statistical tests in univariate statistics, such as the $t$-test \cite{student1908}, the Mann-Whitney test \cite{mann1947}, and the Smola and Scholk\"{o}pf test \cite{smola1998},  are typically utilized when comparing two-sample distributions. 
 In recent years, there has been an increasing interest in multivariate data testing, which plays an important role in real life.
In order to determine whether two multivariate samples drawn from the population distributions $F$ and $G$, respectively, have the same underlying distribution, we formulate such a problem as a two-sample homogeneity test, that is, $$H_{0}: F = G \text{ vs } H_{1}: F \neq G.$$
Some traditional parametric methods, such as multivariate analysis of variance (MANOVA), often rely on the assumption of normality and may not fully capture the complexity of real-world data. This limitation of the normality assumption can be overcome by applying non-parametric tests, such as the Cram\'{e}r test \cite{bar2004} and the energy distance test \cite{szekely2004}, which do not require distribution assumptions to be assigned to the data.
 The biggest challenge with these non-parametric tests is how to extend them to multivariate scenarios. Data depth is a technique for assessing the centrality of multivariate data points in a given data cloud. Let $D(x; F)$ be the data depth function, then the point $x$ in a distribution $F(\cdot)$ from $R^{d}$ can be mapped into a depth value in the interval $[0, 1]$. It has two major advantages: first, it does not require the assumption of data distribution, so it can be used as a non-parametric test; second, the data depth function can be used as a center-outward ranking by mapping multivariate data  to a one-dimensional function unaffected by data dimensionality, and thus facilitates theoretical analysis such as outlier detection and two-sample tests. 
In addition to these advantages, the data depth enjoys four basic properties: affine invariance, central maximality, monotonicity with respect to the deepest point, and vanishing at infinity   \cite{zuo2000}.  
   Some commonly used depth functions, such as  Mahalanobis depth \cite{Liu1993,zuo2006}, Halfspace depth \cite{Dai2023,tukey1975},
   Spatial depth \cite{brown1983, gower1974,zhang2024}, 
   and Projection depth \cite{zuo2000,zuo2003projection}, etc. can be implemented using the \textit{R} package \textit{ddalpha}.

Our research origins from the depth-based quality index $Q(F,G)$, introduced by Liu and Singh \cite{Liu1993},  which measures the relative ``outlyingness'' of the reference distribution $F$ compared to the comparison distribution $G$ and is defined as 
\begin{eqnarray}
Q(F,G):=P\{D(X,F)\leq D(Y,F)| X\sim F,Y\sim G\}=\iint I[D(x; F)\leq D(y; F)]dF(x)dG(y), \label{1}
\end{eqnarray} 
where  $I(\cdot)$ is an indicator function that takes 1 if true and 0 otherwise.
 When $F$ and $G$ are unknown, we denote the empirical distributions of $F$ and $G$ as $F_{m}$ and $G_{n}$, respectively, assuming 
 two
 independent samples are $\{X\}_{i=1}^m$ and $\{Y\}_{j=1}^n$. The $Q$ index can be estimated by the $Q$ statistic
 \begin{eqnarray}
Q(F_{m},G_{n}):=\frac{1}{mn}\sum_{i=1}^m\sum_{j=1}^nI[D(X_i; F_{m})\leq D(Y_j; F_{m})]=\iint I[D(x; F_{m})\leq D(y; F_{m})]dF_{m}(x)dG_{n}(y).\label{2}
\end{eqnarray}

The $Q(F_{m},G_{n})$ and $Q(G_{n},F_{m})$ statistics can make comparisons between samples since it is well known that they are effective approximates of the quality index $Q(F,G)$. However, it is not direct to derive the approximate error. Hence, a deep study of the asymptotic expansion of $Q(F_{m},G_{n})$ and $Q(G_{n},F_{m})$ are necessary. 
So far,  Zuo and He \cite{zuo2006} have made significant contributions to the study of depth-based statistics by showing the asymptotic expansion of  the $Q$ statistic under the null hypothesis. These theoretical insights establish a solid background and foundation for further research.  Shi et al. \cite{shi2023} provided some variants based on   the results of Zuo and He \cite{zuo2006},  including the maximum value of the  $Q$ statistics and the weighted sum of the $Q$ statistics. Recently, Chen et al. \cite{chen2023multivariate}   proposed further variants of the minimum value, the sum, and the product of the $Q$ statistics in order to improve the power. These variants show that $$Q(F_{m},G_{n})+Q(G_{n}, F_{m})-1=O_p(m^{-1})$$ assuming $m/n$  tends to a positive constant for balanced samples. The rate above has been confirmed by Gnettner et al. \cite[Theorem 2.6]{gnettner2023} for the special Halfspace depth  in the case of one-dimensional data. Combining the existing results:
$$Q(F_{m},G_{n})-1/2=O_p(m^{-1/2})\quad\text{and}\quad Q(G_{n}, F_{m})-1/2=O_p(m^{-1/2}),$$   we may obtain that the second term of  $Q(F_{m},G_{n})-1/2$ and $Q(G_{n}, F_{m})-1/2$ is $O_p(m^{-1})$ since the summation of them is $O_p(m^{-1})$. However, for the second term
 Zuo and He \cite{zuo2006} and Gnettner et al.\cite{gnettner2023}  
 show that the optimal rate of the second term for the general depth  is slower than the target one $m^{-1}$.

 In this paper, we shall answer the following questions:
\begin{itemize}
    \item  What are reasonable conditions  for obtaining the target rate  $O_p(m^{-1})$  of the second term in the expansion of the $Q$ statistic based on general depth functions?
    \item      How can current proof methods be extended to obtain a higher-order expansion?
\item How to study the expansion for the available variants?

\end{itemize}
 To answer these questions, we employ the classical Hoeffding decomposition \cite{lee2019}, which has not been used in depth-based $Q$ statistics but is very popular in the study of $U$ statistics. Following the conventional conditions of the Hoeffding decomposition, our new conditions can be not strong enough so that the $Q$ statistic can be decomposed into two uncorrelated main terms and one higher-order term. Moreover, to obtain the rate of  the higher-order term, we extend the Cox-Reid expansion method \citep{cox1987approximations} derived for the approximation of the non-central Chi-squared distribution. Our new expansion results pave the way for the theoretical study of further variants to improve the power.

In Section 2, we introduce the Hoeffding decomposition and the Cox-Reid expansion. Then, we list  the conditions for the theoretical development, revisit  existing methods of proof for  $Q$ statistics,  provide a new method of proof for obtaining an accurate high-order expansion,  provide some examples  to verify our assumptions, and show some variants for increasing power.  
In Section 3, we use two real datasets to illustrate potential applications.  
Finally, we conclude with a brief discussion in Section 4. All lemmas and proofs are placed in the Appendix.
 

\section{Asymptotic expansion of the Q statistic}
For the sake of convenience, we denote:  
\begin{align}
 I(x,y,F)&=I\Big(D(x;F)\leq D(y;F)\Big), ~~\forall~~x,~y~~\in~ ~R^{d},\label{I1}\\
 I(x,y,F_{m},F)&=I(x,y,F_{m})-I(x,y,F), ~~\forall~~x,~y~~ 
 \in~ ~R^{d},\label{I2}\\
 \eta_{l}(Y;F_{m},F)&=E_{X}\Big[\Big(D(Y;F_{m})-D(X;F_{m})\Big)^{l}\Big|D(X;F)=D(Y;F)\Big], l=1,2. \label{I3}  
\end{align} 
\subsection{Hoeffding decomposition}
 
The Hoeffding decomposition is an important tool for studying the properties of $U$-statistic and can be used to decompose the statistic into uncorrelated terms and a higher-order error term.  We apply the Hoeffding decomposition to the study of $Q$ statistics for the first time. The Hoeffding decomposition of $Q$ statistic is:
\begin{eqnarray}
Q(F_{m},G_{n})-Q(F,G)=L_{mn1}+L_{mn2}+R_{mn},\label{pp1}
\end{eqnarray}
where
\begin{eqnarray*}
L_{mn1}&=&\iint I(x,y,F_{m})dF(x)dG_{n}(y)-Q(F,G),\nonumber\\
L_{mn2}&=&\iint  I(x,y,F_{m})dG(y)dF_{m}(x)-Q(F,G),\nonumber\\
R_{mn}
&=& Q(F_{m},G_{n})-Q(F,G)-L_{mn1}-L_{mn2}.
\end{eqnarray*}
Here $L_{mn1}$ and $L_{mn2}$ are the two main terms, which are conditionally independent given some condition related to $F_{m}$, and $R_{mn}$ is the higher order term. We compare the Hoeffding decomposition with the decomposition by 
Zuo and He \cite{zuo2006}:
\begin{eqnarray}
Q(F_{m},G_{n})-Q(F,G)
=S_{mn1} +S_{mn2}+\text{ZH}_{mn},\label{ff11}
\end{eqnarray}
where
\begin{eqnarray*}
 S_{mn1}&=&\iint I(x,y,F)dF(x)d(G_{n}-G)(y),\\
  S_{mn2}&=&\iint I(x,y,F)d(F_{m}-F)(x)dG(y),\\
\text{ZH}_{mn}&=&Q(F_{m},G_{n})-Q(F,G)-S_{mn1}-S_{mn2}.\\
\end{eqnarray*}
It can be seen that the main terms $S_{mn1}$ and $S_{mn2}$ are independent. Since it is stronger than the conditional independence of main terms in Hoeffding decomposition, the rate of remaining term $\text{ZH}_{mn}$ would be slower than the one of $R_{mn}$. Actually, Zuo and He \cite{zuo2006} obtained that 
$\text{ZH}_{mn}=o_{p}(m^{-1/2}).$   Meanwhile, we also compare the Hoeffding decomposition with the decomposition by Gnettner et al. \cite{gnettner2023} as follow:
\begin{eqnarray}
Q(F_{m},G_{n})-Q(F,G)
=SS_{mn}+\text{GKN}_{mn},\label{dd11}
\end{eqnarray}
where
\begin{eqnarray*}
 SS_{mn}&=&\iint I(x,y,F)dF_{m}(x)dG_{n}(y)-\iint I(x,y,F)dF(x)dG(y),\\
\text{GKN}_{mn}&=&Q(F_{m},G_{n})-Q(F,G)-SS_{mn}.\\
\end{eqnarray*}
The main term is a double sum  and the remainder term has an optimal rate $\text{GKN}_{mn}=O_{p}(m^{-3/4})$.
We note that we can obtain a higher rate for $R_{mn}$, i.e., $O_{p}(m^{-1})$,  since we will apply the Cox-Reid expansion to the conditional probabilities to obtain accurate an accurate expansion for $L_{mn1}$ and $L_{mn2}$. 

\subsection{Cox-Reid expansion}
Here, we mainly refer to Theorem 1 of Cox and Reid \cite{cox1987approximations}, which we list below as a lemma:
\begin{lemma}[Cox-Reid expansion]\label{lem 1}
Let $Y=X_{0}+\epsilon X_{1}$ for some fixed $\epsilon>0$. If the joint density of $X_{0}$ and $X_{1}$ satisfies

(i) $\lim\limits _{x_0 \rightarrow-\infty} \frac{\partial^2}{\partial x_0^2} f\left(x_0, x_1\right)=0$,

(ii) $\sup\limits _{x_0} \frac{\partial^2}{\partial x_0^2}\left|f\left(x_0, x_1\right)\right|<K g\left(x_1\right)$, where $\int x_1^3 g\left(x_1\right) d x_1<\infty$,

(iii) $E\left(X_1^2 \mid X_0\right)<\infty$ a.e.,

then
$$
F_Y(y)=F_{X_0}(y)-\epsilon f_{X_0}(y) E\left(X_1 \mid X_0=y\right)+\frac{1}{2} \epsilon^2 \frac{\partial}{\partial y}\left\{f_{X_0}(y) E\left(X_1^2 \mid X_0=y\right)\right\}+O\left(\epsilon^3\right),
$$
where $F_{X}(\cdot)$ is the distribution function of $X$ and $f_{X}(\cdot)$ is the derivative of $F_{X}(\cdot)$.
\end{lemma}
The Cox-Reid expansion provides an accurate approximation of the distribution function of $Y$ in terms of $X_0$ when $\epsilon$ is small. We can apply it to approximate the conditional probabilities for the first integral  as shown in \eqref{2}.  
Under some regular condition   adjusted in the form of data depth (see A5 below), the conditional probability can be expanded as follows: 
\begin{eqnarray}
\int I(x,y,F_{m})dF(x)=F_{D(X;F)}\Big(D(Y;F)\Big)+f_{D(X;F)}\Big(D(Y;F)\Big)\eta_{1}(Y;F_{m},F)+O_{p}(m^{-1}).  \label{z1}
\end{eqnarray}

\subsection{Assumptions}
To show the complete proof, we begin by listing the existing assumptions and add our new reasonable assumptions. 
\begin{itemize}
\item[A1.]    $P(y_{1}\leq D(Y;F)\leq y_{2})\leq C|y_{2}-y_{1}|$ for some positive constants $C$, any $y_{1},y_{2}\in [0,1]$.

\item[A2.]  $\sup\limits_{x\in R^{d}}|D(x;F_{m})-D(x;F)|=o(1)$, almost surely, as $m\rightarrow \infty$.

\item[A3.]    $E\Big(\sup\limits_{x\in R^{d}}|D(x;F_{m})-D(x;F)|^\alpha\Big)=O(m^{-\alpha/2})$ for some $\alpha>0$ (in particular $\alpha=1$).

\item[A4.]    $E(\Sigma_{i}p_{iX}(F_{m})p_{iY}(F_{m}))=o(\rho_{m})$ (in particular $\rho_{m}=m^{-1/2}$) if there exists $c_{i}$ such that $p_{iX}(F_{m})>0$ and $p_{iY}(F_{m}))>0$ for $p_{iZ}(F_{m}))=:P(D(Z;F_{m})=c_{i}|F_{m}), i=1,2,\ldots.$

\item[A5.]    $\eta_{1}(X_{i};F_{m},F)\perp \eta_{1}(X_{k};F_{m},F)| \Gamma_{m}$, $\eta_{1}(Y_{j};F_{m},F)\perp \eta_{1}(Y_{\ell};F_{m},F)| \Gamma_{m}$,\\
 $E\Big[f_{D(Y;F)}\Big(D(X;F)\Big)\eta_{1}(X;F_{m},F)\Big|\Gamma_{m}\Big]=E\Big[f_{D(X;F)}\Big(D(Y;F)\Big)\eta_{1}(Y;F_{m},F)\Big|
\Gamma_{m}\Big]=0$,  \\
$\sup\limits_{\xi}\Big\{\frac{\partial}{\partial(D(X;F))}\Big[f_{D(Y;F)}\Big(D(X;F)\Big)\eta_{2}(X;F_{m},F)\Big]\Big|_{D(X;F)=\xi}\Big\}=O_{p}(m^{-1}),$\\
 $\sup\limits_{\xi}\Big\{\frac{\partial}{\partial(D(Y;F))}\Big[f_{D(X;F)}\Big(D(Y;F)\Big)\eta_{2}(Y;F_{m},F)\Big]\Big|_{D(Y;F)=\xi}\Big\}=O_{p}(m^{-1})$,\\
  where $\Gamma_{m}$ is a condition related to $F_{m}$ and the data depth, and  may be different for different depth functions, $i\neq k$, and $j\neq\ell$. The symbol $\perp$ denotes independence.

\item[A6.]  $E\Big[I(X_{i},Y_{j},F_{m},F)\Big|\Gamma_{m}\Big]=0$ and   
$I(X_{i},Y_{j},F_{m},F)\perp I(X_{k},Y_{l},F_{m},F)| \Gamma_{m}$ for $i\neq k$ and $j\neq l$.
\end{itemize}
\textbf{Remark 1.} In the pioneering study of   Zuo and He \cite{zuo2006},  A1-A4  serve as a basis for asymptotic expansion. Since we need to compare two depth values as shown in 
\eqref{I1} and its difference from the sample version reflected in \eqref{I2}, we may bound the difference of two indicator functions as follows:
 $$|I(x, y, F_m, F)| \leq I\Big(|D(x ; F)-D(y ; F)| \leq 2 \sup _{x \in \mathbb{R}^d}|D(x ; F_m)-D(x ; F)|\Big).$$
   A1   provides an upper bound of its conditional expectation by using the Lipschitz continuity, i.e.,
$$E_Y\left[I\Big(|D(x;F)-D(Y;F)|\leq 2 \sup\limits_{x\in R^{d}}|D(x;F_{m})-D(x;F)|\Big)\right]\leq 4C\sup\limits_{x\in R^{d}}|D(x;F_{m})-D(x;F)|.$$
A2 assumes that the upper bound above is close to zero but without convergence rate.  
A3 further gives a rate on its expectation, which is reasonable according to the standard $\sqrt{m}$ rate of convergence of $F_m$. Combining  A1 with A2 and A3  for a special value of $\alpha=1$ leads to $$\iint I(x,y,F_{m},F)d(F_{m}-F)(x)dG(y)=o_{p}(m^{-\frac{1}{2}}) ~\text{and}~ \iint I(x,y,F_{m},F)dF_{m}(x)d(G_{n}-G)(y)=O_{p}(m^{-\frac{3}{4}}).$$
 The upper bound of the expectation of \eqref{I2} can be estimated by A4 as follows:  
\begin{eqnarray}
\iint I(x, y, F_m, F) d F(x) d G(y)
=\frac{1}{2}E(\Sigma_{i}p_{iX}(F_{m})p_{iY}(F_{m}))=o(\rho_{m}),\nonumber 
\end{eqnarray}
which can be zero for continuous depths.
Finally $\text{ZH}_{mn}=O_p(m^{-3/4})+o_{p}(m^{-1/2})+o(\rho_{m})$, which can be improved to be $O_p(m^{-3/4})$ if we assume $\alpha=2$ and $\rho_{m}=m^{-3/4}$; see more details in Theorem \ref{The 1}.


\noindent\textbf{Remark 2.}
  A5 provides similar conditions in the original Cox-Reid expansion. In particular, it can also lead to $$E\Big[f_{D(Y;F)}\Big(D(X;F)\Big)\eta_{1}(X;F_{m},F)\Big]=E\Big[f_{D(X;F)}\Big(D(Y;F)\Big)\eta_{1}(Y;F_{m},F)\Big]=0.$$    A6 is related to Hoeffding decomposition for accurate expansion better than A1 and implies that $$EI(X_{i},Y_{j},F_{m},F)=E\Big(E_{Y}I(X_{i},Y,F_{m},F)\Big)=E\Big(E_{X}I(X,Y_{j},F_{m},F)\Big)=E\Big(E_{X,Y}I(X,Y,F_{m},F)\Big)=0.$$
We  will give further detailed validations  in \S \ref{2.5}.

\subsection{Asymptotic expansions}
Based on the above assumptions, we   now present the main results as follows.

\begin{theorem}\label{The 1}
Under the null hypothesis $F=G$, A1, A3 with $\alpha=2$, and A4, we have  $\text{ZH}_{mn}=O_{p}(m^{-3/4})+o(\rho_{m}).$
\end{theorem}
The proof of Theorem \ref{The 1} is provided in Appendix 6.2. By a standard central limit theorem, we have
$$\sqrt{n}S_{mn1}\stackrel{d}\longrightarrow N(0, {1}/{12})~\text{and}~
\sqrt{m}S_{mn2}\stackrel{d}\longrightarrow N(0, {1}/{12}),\quad\text{as}~m, n\rightarrow\infty,$$
where $S_{mn1}$ and $S_{mn2}$ are defined in \eqref{ff11}. Therefore, we have the following corollary.
\begin{corollary}\label{cor 1}
Under the conditions of Theorem \ref{The 1} and $\rho_{m}=m^{-1/2}$, as $m, n\rightarrow\infty$,
\begin{eqnarray}
\Big[\frac{1}{12}\left(\frac{1}{m}+\frac{1}{n}\right)\Big]^{-\frac{1}{2}}\Big[Q(F_{m},G_{n})-\frac{1}{2}\Big]\stackrel{d}\longrightarrow N(0,1).\nonumber
\end{eqnarray}
\end{corollary}
\noindent\textbf{Remark 3.} 
  Theorem 1 of Zuo and He \citep{zuo2006} is based on  A1, A2, A3 with $\alpha=1$, and A4, and hence $\text{ZH}_{mn}=o_{p}(m^{-1/2})+o(\rho_{m}).$ 
  In fact,  if we choose $\alpha=2$ in A3, then the rate $o_{p}(m^{-1/2})$ can be improved to   $O_{p}(m^{-3/4})$, which coincides with the best rate in Gnettner et al. \cite{gnettner2023}. Next, we will use the Hoeffding decomposition and Cox-Reid expansion 
to obtain the faster rate $O_{p}(m^{-1})$.


\begin{theorem}\label{The 2}
Under the null hypothesis $F=G$, A3 with $\alpha=2$, A4-A6, we have
\begin{align}
&Q(F_{m},G_{n})-\frac{1}{2}=\frac{1}{n}\sum\limits_{j=1}^{n}\Big[F_{D(X;F)}\Big(D(Y_{j};F)\Big)-\frac{1}{2}\Big]+\frac{1}{m}\sum\limits_{i=1}^{m}\Big[\frac{1}{2}-F_{D(Y;F)}\Big(D(X_{i};F)\Big)\Big]+R_{mn},\label{ww1}\\
&Q(G_{n},F_{m})-\frac{1}{2}=\frac{1}{n}\sum\limits_{j=1}^{n}\Big[\frac{1}{2}-F_{D(X;F)}\Big(D(Y_{j};F)\Big)\Big]+\frac{1}{m}\sum\limits_{i=1}^{m}\Big[F_{D(Y;F)}\Big(D(X_{i};F)\Big)-\frac{1}{2}\Big]+\tilde{R}_{mn},\label{w1}
\end{align}
where $R_{mn}=\tilde{R}_{mn}=O_{p}(m^{-1})+o(\rho_{m})$.
\end{theorem}

\noindent\textbf{Remark 4.} Theorem \ref{The 2} 
shows that the normalized $Q$ statisticss have the same asymptotic distribution:  $$\Big[\frac{1}{12}\left(\frac{1}{m}+\frac{1}{n}\right)\Big]^{-\frac{1}{2}}\Big[Q(F_{m},G_{n})-\frac{1}{2}\Big]\stackrel{d}\longrightarrow N(0,1),\quad \Big[\frac{1}{12}\left(\frac{1}{m}+\frac{1}{n}\right)\Big]^{-\frac{1}{2}}\Big[Q(G_{n},F_{m})-\frac{1}{2}\Big]\stackrel{d}\longrightarrow N(0,1),$$ which suggests that both distributions $F$ and $G$ are significantly different at  level 0.05   if the absolute value of the normalized $Q$ statistic is greater than 1.96. 
In addition,   Theorem \ref{The 2} gives insights for studying the asymptotic distributions of so-called sum and product of Q statistics
\citep{chen2023multivariate,chen2024}. For continuous depths or $\rho_m=O(m^{-1})$,  we obtain that $$Q(F_{m},G_{n})+Q(G_{n},F_{m})-1=O_{p}(m^{-1})\quad\text{and}\quad\Big[Q(F_{m},G_{n})-1/2\Big]\Big[Q(G_{n},F_{m})-1/2\Big]=O_{p}(m^{-1}).$$

 

\subsection{Verifications of assumptions}\label{2.5}
To illustrate the wide application of our results, we verify the conditions for these five depth functions: Euclidean depth, Mahalanobis depth, Halfspace depth,  Projection depth, and Spatial depth.  
 Since the Euclidean depth is a special case of the Mahalanobis depth, Zuo and He \cite{zuo2006} have already discussed A1-A4 for the Mahalanobis depth, the half-space depth, and the projection depth,  we only need to check the rest of the assumptions A5-A6 for these three depths, and all of the assumptions A1-A6 for the Spatial depth.

\begin{exmp}
[Euclidean depth \cite{kosiorowski2014depthproc}]  The univariate Euclidean depth of any point $x\in R$ in a one-dimensional distribution $F$ is
defined as
\begin{eqnarray}
ED(x;F)=\frac{1}{1+(x-\mu_{F})^2},\label{mm0}
\end{eqnarray}
where $\mu_{F}$ is the mean of the distribution $F$. 
\end{exmp}
The empirical version of $ED(X;F)$ is $ED(X;F_{m})=1/(1+(X-\Bar{X})^2)$, where $\Bar{X}=\sum_{i=1}^{m}X_{i}/m$. Since $X_{i}$ and $X_{k}$ are independent for $i\neq k$, it follows that $\eta_{1}(X_{i};F_{m},F)$ is independent of $\eta_{1}(X_{k};F_{m},F)$ when given $\Gamma_{m}$. The same is true for the rest of the examples.
While Euclidean depth is simple for one-dimensional data, we can look closely and find connections to other depths. The first step is   to simplify the calculation of $f_{D(Y;F)}\Big(D(X;F)\Big)$ involved in A5. To do this, we convert the Euclidean depth to the Euclidean distance $d(x;F)=|x-\mu_{F}|$ as follows.
  $$ 
F_{ED(X;F)}\Big(ED(Y;F)\Big)=P\Big(\frac{1}{1+(X-\mu_{F})^2}\leq \frac{1}{1+(Y-\mu_{F})^2}\Big)=1-F_{d(X;F)}\Big(d(Y;F)\Big).$$ 
We assume that the density function $f_{d(X;F)}(\cdot)$ is bounded from above.
Since $ED(X;F_{m})$ does not depend on the mean, we further assume that $\mu_{F}=0$ without loss of generality.  
We consider the condition   $\Gamma_m$  to be $\bar X$, reflecting the information in $\mu_F$.
Since $f_{d(X;F)}\Big(d(Y;F)\Big)$ is bounded from above, $d(Y;F)$ is an even function of $Y$, and   $\eta_{1}(Y;F_{m},F)|\Gamma_m$ is equal to zero or an odd function of $Y$, we have
 $$E_Y\Big[f_{d(X;F)}\Big(d(Y;F)\Big)\eta_{1}(Y;F_{m},F)\Big|\Gamma_{m}\Big]=0.$$
Since $\eta_{2}(Y;F_{m},F)=O_p(m^{-1})$, we use Cauchy's theorem to bound the partial derivatives in terms of integrals, where the integrating function is $O_p(m^{-1})$, leading to 
$$
\sup\limits_{\xi}\Big\{\frac{\partial}{\partial(d(Y;F))}\Big[f_{d(X;F)}\Big(d(Y;F)\Big)\eta_{2}(Y;F_{m},F)\Big]\Big|_{d(Y;F)=\xi}\Big\}= O_{p}(m^{-1}).$$ 
Therefore,  A5 holds. To verify A6, we first consider  $i\neq k \in \{1,\cdots,m\}$ and $j \neq l\in \{1,\cdots,n\}$ and have that 
$$ E\Big\{E\Big[I(X_{i},Y_{j},F_{m},F)I(X_{k},Y_{l},F_{m},F)\Big|\Gamma_{m}\Big]\Big\}\nonumber\\
=E\Big\{E\Big[I(X_{i},Y_{j},F_{m},F)\Big|\Gamma_{m}\Big]E\Big[I(X_{k},Y_{l},F_{m},F)\Big|\Gamma_{m}\Big]\Big\}. $$
 Moreover, noting that $ED(X_{i};F_{m})$ and $ED(Y_{j};F_{m})$ are conditionally independent and have the same distribution given $\Gamma_m$, we have  
 $$E\Big[I(X_{i},Y_{j},F_{m},F)\Big|\Gamma_{m}\Big]=P\Big(\Big(ED(X_{i};F_{m})\leq ED(Y_{j};F_{m})\Big)\Big|\Gamma_{m}\Big)-P\Big(ED(X_{i};F)\leq ED(Y_{j};F)\Big)=0.$$
 
\begin{exmp}[Mahalanobis depth   \cite{zuo2006}]   For any point $x\in R^{d}$, the Mahalanobis depth of   point $x$ is defined follows:
\begin{eqnarray}
MD(x;F)=\frac{1}{1+(x-\mu_{F})'\Sigma^{-1}_{F}(x-\mu_{F})},\label{m1}
\end{eqnarray} 
where $F$ is a d-dimensional distribution and $\mu_{F}$ and $\Sigma_{F}$ are the mean and covariance of $F$, respectively.

\end{exmp}
The empirical version  of $MD(X;F)$ is $MD(X;F_{m})=\frac{1}{1+(X-\Bar{X})'\Sigma_{F_{m}}^{-1}(X-\Bar{X})}$, where $\Bar{X}$ and $\Sigma_{F_{m}}$ are the sample mean and sample covariance, respectively. Similarly, we convert the Mahalanobis depth to   the Mahalanobis distance $d(x;F)=\|\Sigma^{-\frac{1}{2}}_{F}(x-\mu_{F})\|_2$  as follows.
 $$ 
F_{MD(X;F)}\Big(MD(Y;F)\Big)=1-F_{d(X;F) }\Big(d(Y;F)\Big).  $$
We also assume that the density function $f_{d(X;F)}(\cdot)$ is bounded from above and that $\mu_{F}=0$. Here, we consider the condition $\Gamma_m$ to be $\bar X$ and $\Sigma_{F_{m}}$. Again, since  $\eta_{1}(Y;F_{m},F)|\Gamma_m$ is equal to zero or an odd function of $Y$, we have
\begin{eqnarray}
E_Y\Big[f_{d(X;F)}\Big(d(Y;F)\Big)\eta_{1}(Y;F_{m},F)\Big|\Gamma_{m}\Big]=0
.\nonumber
\end{eqnarray}
Notice that $\eta_{2}(Y;F_{m},F)=O_p(m^{-1})$, and the rest of the argument follows in the same way. 

\begin{exmp}[Halfspace depth \cite{Dai2023,zuo2006}]
For any point $x\in R^{d}$, the Halfspace depth, also known for the Tukey depth, is defined as 
\begin{eqnarray}
HD(x;F)=\inf\{P(H_{x}):~H_{x}~ is~ a~ closed~ half~ space~ containing~ x \},\label{m0}
\end{eqnarray}
where $P$ is the probability measure corresponding to $F$. 
\end{exmp} 
If we replace the probability $P$ by the empirical probability, we obtain the sample version of $HD(x, F)$, denoted $HD(x, F_m).$
  In light of  Zuo and He \cite{zuo2006,Dai2023}, the same argument is $\eta_{2}(Y;F_{m},F)=O_p(m^{-1})$,    but the different things are that $f_{HD(X;F)}(HD(Y;F))$ is a constant and $\Gamma_m=F_m$ leading to  $$E\Big[\Big(HD(Y; F_m)-HD(X; F_m)\Big)\Big|F_m\Big]=0.$$
Therefore, we have
$$E\Big[f_{HD(X;F)}\Big(HD(Y;F)\Big)\eta_{1}(Y;F_{m},F)\Big|\Gamma_{m}\Big]=0.$$
The remain arguments are the same. 
\begin{exmp}[Projection depth  \cite{zuo2003projection}] For any point $x\in R^{d}$, the Projection depth is defined as 
$$
P D(x ; F)=1 /(1+O(x ; F)),
$$
where  
$
O(x; F)=\sup _{u \in S^{d-1}}\left|u^{\prime} x-\mu\left(F_u\right)\right| / \sigma\left(F_u\right),
$, $S^{d-1}=\{u:\|u\|=1\}, \mu(F)$ and $\sigma(F)$ are   location and scale of $F$, respectively,  and $u^{\prime} X \sim F_u$ with $X \sim F$.  
\end{exmp}

The empirical version  of $PD(X;F)$ is $PD(X;F_{m})=\frac{1}{1+\sup _{u \in S^{d-1}}\left|u^{\prime} X-\mu\left(F_{mu}\right)\right| / \sigma\left(F_{mu}\right)}$, where $F_{mu}$ is the empirical distribution of $u^{\prime} X_1,\cdots,u^{\prime} X_m$. 
 As shown in Zuo and He \cite{zuo2006}, unlike the Halfspace depth, we will assume that
$\mu(a X+b)=a \mu(X)+b$ and $\sigma(a X+b)=|a| \sigma(X)$ for any scalars $a, b \in  {R}^1$.
Let $\Gamma_{m}$ be $\mu\left(F_{mu}\right)$ and
 $\sigma\left(F_{mu}\right)$, the rest of the verification is the same.
 
\begin{exmp}[Spatial depth \cite{serfling2002,zhang2024}]
For any point $x\in R^{d}$, the Spatial depth is defined as 
$$SD(x;F)=1-\Big\|E_X\frac{(x-X)}{\|x-X\|}\Big\|,\quad X\sim F.$$  
\end{exmp} 
The empirical version  of $SD(X;F)$ is $SD(X;F_{m})=1-\Big\|\frac{1}{m}\sum\limits_{i=1}^{m}\frac{(X-X_{i})}{\|X-X_{i}\|}\Big\|$, which does not depend  on the mean of $F$.   
We assume that $X$ is continuous and its density is symmetric around mean zero, the density of $SD(X;F)$ denoted as $f_{SD(X;F)}(\cdot)$ is bounded from above, and $E\sup\limits_{x\in R^{d}}\Big\|\sum\limits_{i=1}^{m}\Big(\frac{(x-X_{i})}{\|x-X_{i}\|}-E\frac{(x-X_{i})}{\|x-X_{i}\|}\Big)\Big\|^{\alpha}=O(m^{\alpha/2})$.  
Then, we can conclude that $SD(X;F)$ is continuous, so A1 holds  and A4 holds at $\rho_m=0$. 
By using the inequality
$$ 
\sup\limits_{x\in R^{d}}|SD(x;F_{m})-SD(x;F)| \leq\frac{1}{m}\sup\limits_{x\in R^{d}}\Big\|\sum\limits_{i=1}^{m}\Big(\frac{(x-X_{i})}{\|x-X_{i}\|}-E\frac{(x-X_{i})}{\|x-X_{i}\|}\Big)\Big\|,$$  
A2 holds by the strong law of large number theory and A3 holds.  
 Similarly, $SD(Y; F)$ is an even function of $Y$ because $X$ and $-X$ are identically distributed  and  $\eta_{1}(Y;F_{m},F)|\Gamma_m$ is equal to zero or an odd function of $Y$, 
$$E_Y\Big[f_{SD(X;F)}\Big(SD(Y;F)\Big)\eta_{1}(Y;F_{m},F)\Big|\Gamma_{m}\Big]=0.$$
Again, A5 follows by $\eta_{2}(Y;F_{m},F)=O_p(m^{-1})$,  setting  $\Gamma_m=X_1, \ldots, X_m$, and redefining   
$\frac{X_{i}-X_{t}}{\|X_{i}-X_{t}\|}$ at $i=t$ to $\frac{Y_{j}-X_{t}}{\|Y_{j}-X_{t}\|}$  for some introduced $Y_{j}$. In addition, A6 holds by  
$$E\Big[I(X_{i},Y_{j},F_{m},F)\Big|\Gamma_{m}\Big]=P\Big(\Big(SD(X_{i};F_{m})\leq SD(Y_{j};F_{m})\Big)\Big|\Gamma_{m}\Big)-P\Big(SD(X_{i};F)\leq SD(Y_{j};F)\Big)
=0.$$

\subsection{Some variants}
 
It is worth noting that $Q(F, G)=Q(G, F)=1/2$ under the null hypothesis $F=G$ but  $Q(F_m, G_n)$ and $Q(G_n, F_m)$ may not be same. Theorem 
\ref{The 2} implies  
$$Q(F_m, G_n)-1/2)=-\Big[Q(G_n, F_m)-1/2\Big]+O_p(m^{-1})+o(\rho_m).$$

To obtain more power under alternative hypothesis, \cite{shi2023}  proposed the maximum statistic $M_{m,n}$ defined as
\begin{eqnarray}\label{mq1}
    M_{m,n}=\Big[\frac{1}{12}\Big(\frac{1}{m}+\frac{1}{n}\Big)\Big]^{-1}\max\Big[\Big(Q(F_{m},G_{n})-1/2\Big)^2, \Big(Q(G_{n},F_{m})-1/2\Big)^2\Big].
\end{eqnarray} 
\cite{chen2023multivariate}  proposed the minimum statistic written as
\begin{eqnarray}\label{mq2}
M^{*}_{m,n}
=\Big[\frac{1}{12}\Big(\frac{1}{m}+\frac{1}{n}\Big)\Big]^{-1/2}\Big
[1/2-\min\Big(Q(F_{m},G_{n}),Q(G_{n},F_{m})\Big)\Big].\end{eqnarray} 
Both variants $M_{m,n}$ and $M^{*}_{m,n}$ are so-called same-attraction functions and $M_{m,n}$ and $(M^{*}_{m,n})^2$ have the same attractor \citep{chen2024}, where  “attractor”
refers to the convergence of the distribution to a particular limit as the sample size increases, while “attraction” denotes the direction in which the terms  move in $\max$ and $\min$, respectively.

We have further expansions of those two variants as follows.
\begin{corollary}\label{cor 2}
Under the assumptions of Theorem \ref{The 2}, we have
\begin{align}
&M_{m,n}=\Big[\frac{1}{12}\Big(\frac{1}{m}+\frac{1}{n}\Big)\Big]^{-1}\Big\{\frac{1}{n}\sum\limits_{j=1}^{n}\Big[F_{D(x;F)}\Big(D(Y_{j};F)\Big)-\frac{1}{2}\Big]+\frac{1}{m}\sum\limits_{i=1}^{m}\Big[\frac{1}{2}-F_{D(y;F)}\Big(D(X_{i};F)\Big)\Big]\Big\}^2+\tilde{R}_{mn1},\nonumber\\
&M^{*}_{m,n}
=\Big[\frac{1}{12}\Big(\frac{1}{m}+\frac{1}{n}\Big)\Big]^{-1/2}\Big|\frac{1}{n}\sum\limits_{j=1}^{n}\Big[F_{D(x;F)}\Big(D(Y_{j};F)\Big)-\frac{1}{2}\Big]+\frac{1}{m}\sum\limits_{i=1}^{m}\Big[\frac{1}{2}-F_{D(y;F)}\Big(D(X_{i};F)\Big)\Big]\Big|+\tilde{R}_{mn2},\nonumber
\end{align}
where $\tilde{R}_{mn1}=O_{p}(m^{-1/2})+o_p(m^{1/2}\rho_{m})$ and $\tilde{R}_{mn2}=O_{p}(m^{-1/2})+o(m^{1/2}\rho_{m})$.
\end{corollary} 
The following corollary further confirms that both variants have the same attractor.
\begin{corollary}\label{cor 3}
Under the assumptions of Theorem \ref{The 2}, we have
\begin{eqnarray}
M_{m,n} \xrightarrow d \chi_{1}^2~~and~~(M^{*}_{m,n})^2\xrightarrow d \chi_{1}^2, \quad\text{as}~m, n\rightarrow\infty.\nonumber
\end{eqnarray}

\end{corollary}
Therefore, we  suggest  that both distributions $F$ and $G$ are significantly different at  level 0.05   if $M_{m,n}>1.96^2$ and $(M^{*}_{m,n})^2>1.96^2$, respectively.

\section{Data Analysis}
\subsection{Cancer image data}
 
Breast cancer is the second leading cause of cancer deaths among women in Canada. Early detection reduces the cost of treatment and provides patients with a favorable prognosis. Traditional methods are labor-intensive and therefore consume significant healthcare resources. Recent studies have complemented traditional methods using automated mammography analysis.
We consider a digital mammogram with potentially tumours from the VinDr-Mammo database \cite{cancerdata}, as shown in Fig. \ref{cancer} (a).    
A fuzzy c-means clustering algorithm based on Euclidean distance was applied to this mammogram to evaluate the tumor detection effect, setting the fuzzy degree to 2 and the number of clusters to 3 \citep{krasnov}. When each point was clustered into a  cluster with the highest probability, three different clusters are shown in Fig. \ref{cancer} (b-d).
  
\begin{center}
\begin{tabular}{ccccc}
\includegraphics[width=0.2\linewidth]{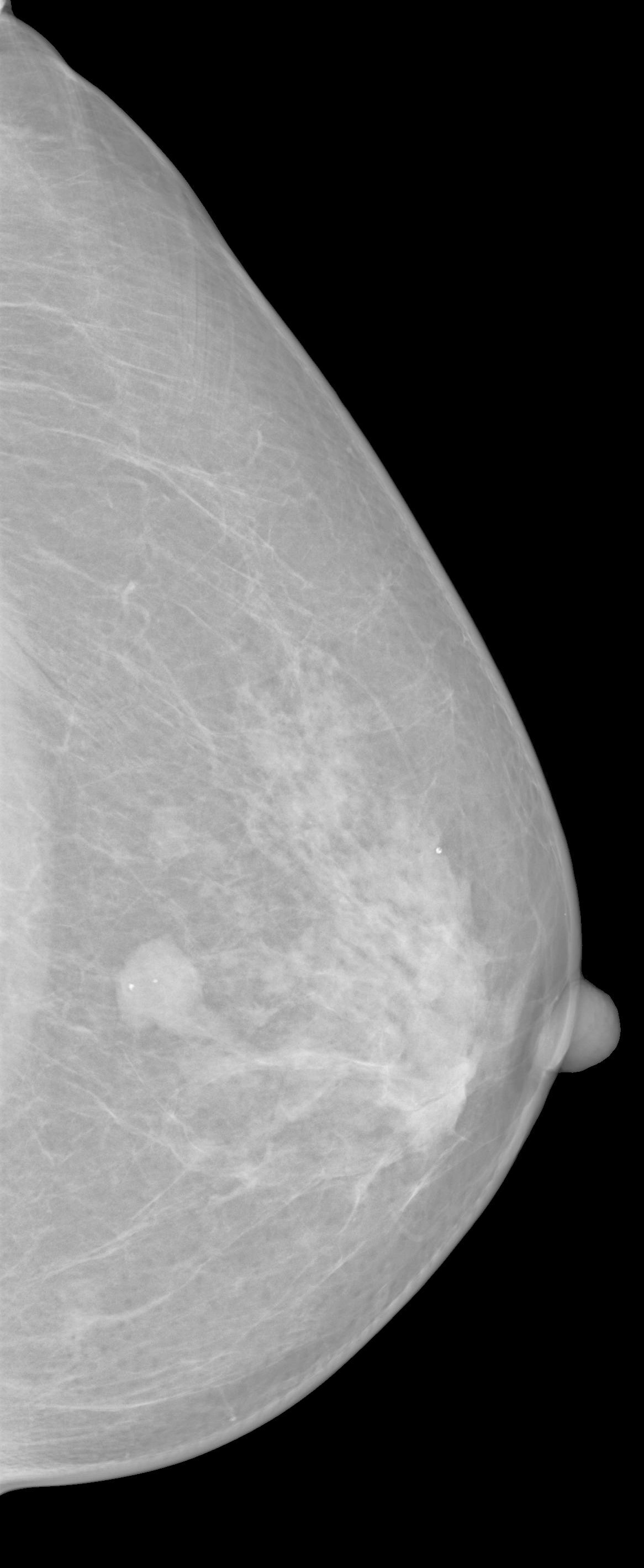} &
\includegraphics[width=0.2\linewidth]{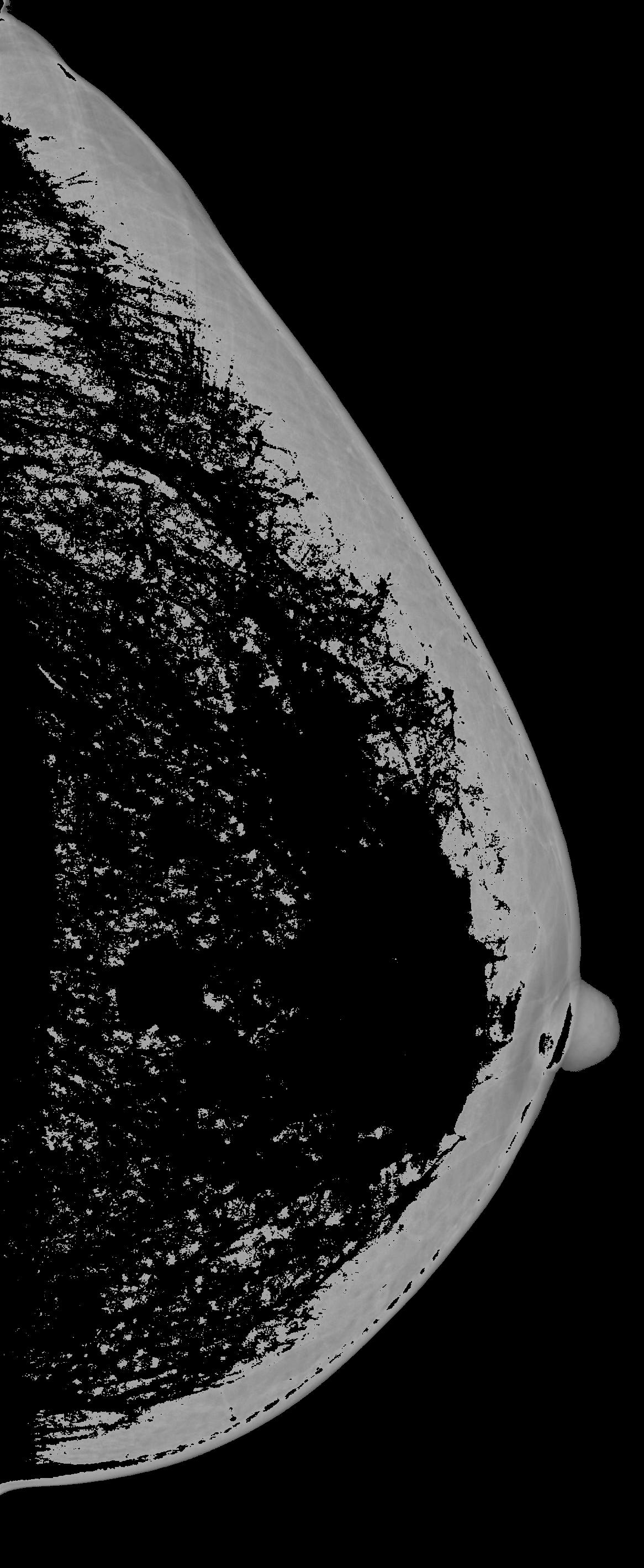} &
\includegraphics[width=0.2\linewidth]{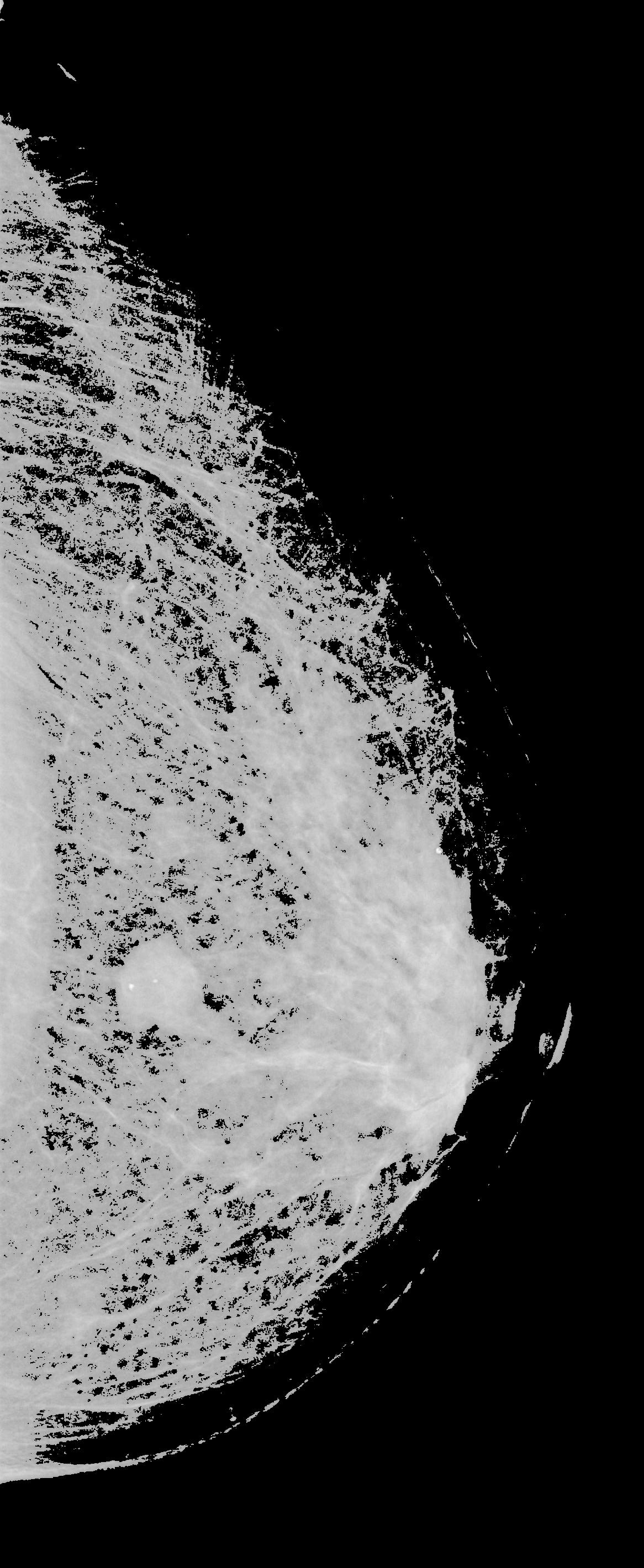} &
\includegraphics[width=0.2\linewidth]{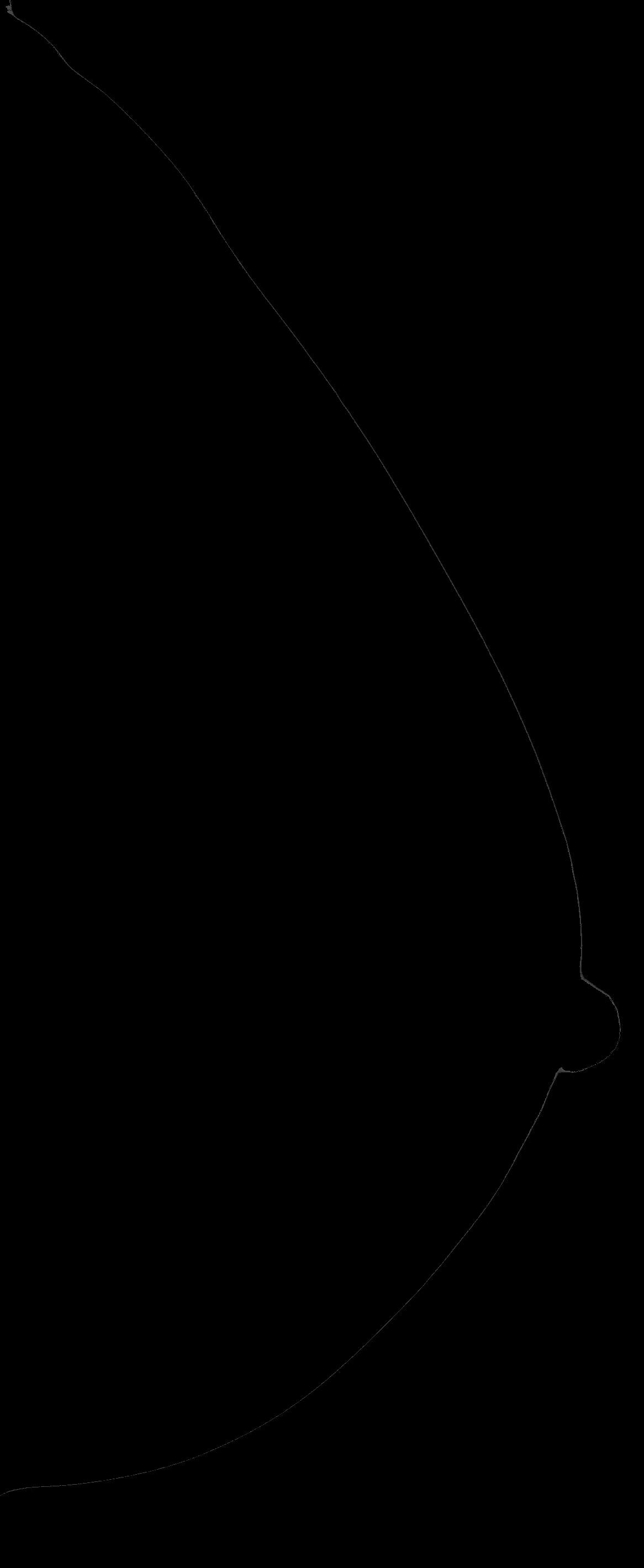} \\
 (a): original image & (b): cluster 1 & (c): cluster 2 & (d): cluster 3.  \\
\end{tabular}
\captionof{figure}{Euclidean-distance-based fuzzy c-means clustering on the breast cancer image. (a): original image; (b): cluster 1; (c): cluster 2; (d): cluster 3}
\label{cancer}
\end{center}

It is well known that the fuzzy c-means clustering algorithm is a soft clustering \citep{fcm} algorithm that allows data samples to belong to multiple clusters. However, the above clustering is based on hard fuzzy clustering, where each data point can only belong to one cluster with the highest probability. 
To confirm whether hard clustering is accurate, we apply a two-sample test to further validate the clustering results as shown in Fig. \ref{cancer}. There are 862532, 1404409 and 26345 data observations in three different clusters. Our aim is to test the difference between cluster 1 and cluster 2 as they are close to each other. We consider the two $Q$ statistics in \eqref{ww1} and \eqref{w1}  denoted as $Q_{m,n}$ and $Q_{m,n}^*$, respectively,  the maximum statistic  $M_{m,n}$ in \eqref{mq1}, and the minimum statistic  $M _{m,n}^*$ in \eqref{mq2}. 
For the Euclidean depth, the Halfspace depth, and the Projection depth, the corresponding four tests are performed. At these depths, all tests have an asymptotic $p$-value of 0, indicating a significant difference between cluster 1 and cluster 2.

\subsection{Skull data}

The brains of humans may change over time. The reason is still puzzling researchers. To examine changes in skull size over time, we consider the  Egyptian skull data  obtained from the \textit{R} package \texttt{HSAUR}, consisting of four measurements of skulls (maximum breaths, basibregmatic heights, basialiveolar length, and nasal heights of the skull). We focus on two epochs, 1850 B.C. and 200 B.C.. Each epoch contains 30 samples. We investigate whether skull size changes over time during interbreeding with immigrants.

\begin{figure} 
\begin{center}
 \includegraphics[width=0.7\textwidth]{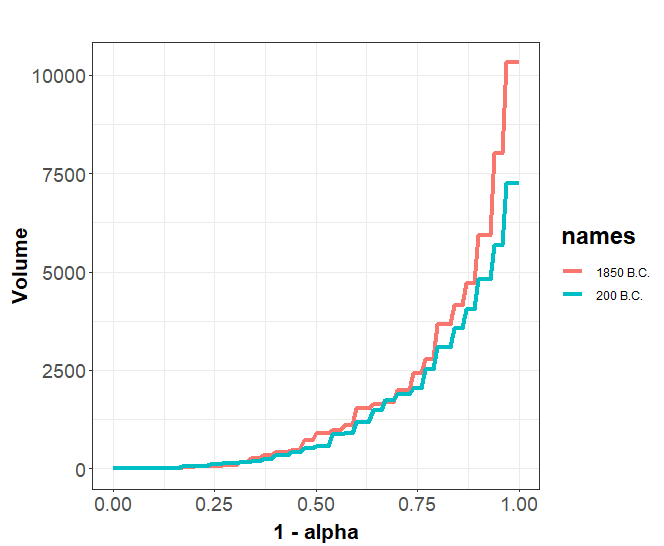}
\caption{Scale curves of skull data for epochs: 1850 B.C. and 200 B.C. under Mahalanobis depth} 
 \label{skull-3-4}
\end{center}
\end{figure}

First, we plot the scale curves \cite{Liu99} to visualize the differences in skull sizes, shown in Figure \ref{skull-3-4}. The scale curve compares the dispersion of two samples.
As observed from the scale curve, these two curves show a small difference. We further apply the four depth-based tests above. The asymptotic $p$-values are shown in  Table \ref{table-skull-3-4}.

\begin{table}[H]
\centering
\begin{tabular}{ |p{3cm}| p{1cm} |p{1cm}| p{1cm}|p{1cm}|} 
\hline 
  & $Q_{m,n}$ & $Q_{m,n}^*$ &  $M_{m,n}$ &   $M_{m,n}^*$ \\
\hline 
 Mahalanobis depth &0.676 & 0.051   &  0.051 & 0.051 \\ 
 Halfspace depth & 0.917 & 0.149 &   0.149 & 0.149 \\ 
 Projection depth & 0.982 & 0.282 &  0.282 & 0.282  \\ 
 Spatial depth&0.996 &0.432&0.432&0.432\\
\hline 
\end{tabular}
\caption{Asymptotic $p$-values for $Q_{m,n}$, $Q_{m,n}^*$, $M_{m,n}$, $M_{m,n}^*$  under different depths.}
\label{table-skull-3-4}
\end{table}
  
We can see that the asymptotic $p$-values for all depths are greater than 0.05, suggesting that there is no strong correlation between skull size and immigrant hybridization in both epochs.

\section{Conclusion and discussion}

By revisiting the existing assumptions of the asymptotic analysis of the $Q$ statistic, we provide two additional plausible assumptions to obtain higher-order expansions based on the Hoeffding decomposition and the Cox-Reid expansion. Our results are applicable to a number of popular depths and converge faster than existing results. In addition, we discuss some variants and show their expansions that  are not known so far. The potential applications are illustrated by analyzing cancer images and skull size variations. Although we only provide expansions to the sum  and product of $Q$ statistics, our expansions provide insights 
for further study of their asymptotic distributions.


\section{Appendix}
\subsection{Main Lemma and Proof of Lemma}
\begin{lemma}[Lemma A.2 \cite{gnettner2023}] \label{lem 2}
Suppose $H_{1}$ and $H_{2}$ are arbitrary distribution functions, for any point $x,y,z \in R^{d}$, it holds
\begin{eqnarray}
|I(x,y,H_{1})-I(x,y,H_{2})|\leq I\Big(|D(x,H_{1})-D(y,H_{1})|\leq 2\sup\limits_{z\in R^{d}}|D(z,H_{1})-D(z,H_{2})|\Big).\nonumber
\end{eqnarray}
\end{lemma}
\textbf{Proof of Lemma \ref{lem 2}.} Observing  $I(x,y,H)=I(D(x,H)\leq D(y,H))$, if the left side equals 0, then the inequality always holds true. Therefore, we only need to consider the following cases: $|I(x,y,H_{1})-I(x,y,H_{2})|=1$ when $D(x,H_{1})\leq D(y,H_{1})$ and $D(x,H_{2})>D(y,H_{2})$ or vice versa.  
\begin{eqnarray}
|D(x,H_{1})-D(y,H_{1})|&=&D(x,H_{1})-D(y,H_{1})\nonumber\\
&\leq& D(x,H_{1})-D(y,H_{1})+D(y,H_{2})-D(x,H_{2})\nonumber\\
&=&\Big(D(y,H_{2})-D(y,H_{1})\Big)+\Big(D(x,H_{1})-D(x,H_{2})\Big)\nonumber\\
&\leq& 2\sup\limits_{z\in R^{d}}|D(z,H_{1})-D(z,H_{2})|.\nonumber
\end{eqnarray}
In a similarly way, we can obtain the proof for other cases. The proof for Lemma \ref{lem 2} is complete.\\

Lemma \ref{lem 2} plays a significant role in our proof. However, Zuo and He \cite{zuo2006} did not provide a proof. The important role of Lemma \ref{lem 2} was also noted by Gnettner et al. \cite{gnettner2023}  and a related proof was provided in their Lemma A.2.

\begin{lemma}\label{lem 3}
Under   A1 and A3 with $\alpha=1$, we have

\begin{eqnarray}
\iint I(x,y,F_{m},F)dF(x)d(G_{n}-G)(y)=O_{p}\Big(m^{-1/4}n^{-1/2}\Big),\label{a0}
\end{eqnarray}
\begin{eqnarray}
\iint I(x,y,F_{m},F)dF_{m}(x)d(G_{n}-G)(y)=O_{p}\Big(m^{-1/4}n^{-1/2}\Big).\label{a00}
\end{eqnarray}
\end{lemma}
\textbf{Proof of Lemma \ref{lem 3}.}
Let $W_{nm1}=\iint I(x,y,F_{m},F)dF(x)d(G_{n}-G)(y)$. By the independence of the sequence, we have
\begin{eqnarray}
EW_{nm1}&=&E\iint I(x,y,F_{m},F)dF(x)d(G_{n}-G)(y)\nonumber\\
&=&E\Big[\frac{1}{n}\sum_{j=1}^{n}\Big(\int I(x,Y_{j},F_{m},F)dF(x)-E(I(X,Y,F_{m},F))\Big)\Big]\nonumber\\
&=&E\Big[\frac{1}{n}\sum_{j=1}^{n}\Big(E_{X}(I(X,Y_{j},F_{m},F))-E(I(X,Y,F_{m},F))\Big)\Big]\nonumber\\
&=&E\Big(E_{X}(I(X,Y_{1},F_{m},F))-E(I(X,Y,F_{m},F))\Big)=0,\nonumber
\end{eqnarray}
and by Cauchy-Schwarz inequality, we get
\begin{eqnarray}
EW^{2}_{nm1}&=&E\Big\{E\Big[\Big(\iint I(x,y,F_{m},F)dF(x)d(G_{n}-G)(y)\Big)^{2}\Big|F_{m}\Big]\Big\}\nonumber\\
&\leq&E\Big\{E\Big[\int\Big(\int I(x,y,F_{m},F)d(G_{n}-G)(y)\Big)^{2}dF(x)\Big|F_{m}\Big]\Big\}\nonumber\\
&=&E\Big\{E_{X}\Big[\Big(\int I(X,y,F_{m},F)d(G_{n}-G)(y)\Big)^{2}\Big|F_{m}\Big]\Big\}\nonumber\\
&=&E\Big\{E_{X}\Big[\Big(\frac{1}{n}\sum_{j=1}^{n}I(X,Y_{j},F_{m},F)-E_{Y}I(X,Y,F_{m},F)\Big)^{2}\Big|F_{m}\Big]\Big\}\nonumber\\
&\leq&\frac{1}{n}E\Big\{E\Big[\Big(I(X,Y_{1},F_{m},F)\Big)^{2}\Big|F_{m}\Big]\Big\}\nonumber\\
&=&\frac{1}{n}E\Big\{E\Big[|I(X,Y_{1},F_{m},F)|\Big|F_{m}\Big]\Big\}.\label{a111}
\end{eqnarray}
By Lemma \ref{lem 2},  A1 and A3 with $\alpha=1$, there is
\begin{eqnarray}
EW^{2}_{nm1}\leq\frac{4C}{n}E\Big(\sup\limits_{x\in R^{d}}|D(x;F_{m})-D(x;F)|\Big)=O\Big(m^{-1/2}n^{-1}\Big).\label{a2}
\end{eqnarray}
By using the Chebyshev's inequality, we can obtain
\begin{eqnarray}
W_{nm1}=O_{p}\Big(m^{-1/4}n^{-1/2}\Big).\label{a3}
\end{eqnarray}
Similarly, let $W_{nm2}=\iint I(x,y,F_{m},F)dF_{m}(x)d(G_{n}-G)(y)$,
\begin{eqnarray}
EW_{nm2}&=&E\iint I(x,y,F_{m},F)dF_{m}(x)d(G_{n}-G)(y)\nonumber\\
&=&E\Big(\frac{1}{m}\sum_{i=1}^{m}\int I(X_{i},y,F_{m},F)d(G_{n}-G)(y)\Big)\nonumber\\
&=&\frac{1}{m}\sum_{i=1}^{m}\Big\{E\Big[\frac{1}{n}\sum_{j=1}^{n}\Big(I(X_{i},Y_{j},F_{m},F)-E_{Y}I(X_{i},Y,F_{m},F)\Big)\Big]\Big\}\nonumber\\
&=&0.\nonumber
\end{eqnarray}
By Cauchy-Schwarz inequality, it has
\begin{eqnarray}
EW^{2}_{nm2}&=&E\Big[\iint I(x,y,F_{m},F)dF_{m}(x)d(G_{n}-G)(y)\Big]^{2}\nonumber\\
&\leq& E\Big[\int\Big(\int I(x,y,F_{m},F)d(G_{n}-G)(y)\Big)^{2}dF_{m}(x)\Big]\nonumber\\
&=&E\Big[\int I(X_{1},y,F_{m},F)d(G_{n}-G)(y)\Big]^{2}\nonumber\\
&=&E\Big[E\Big(\Big(\frac{1}{n}\sum_{j=1}^{n}I(X_{1},Y_{j},F_{m},F)-E_{Y}I(X_{1},Y,F_{m},F)\Big)^{2}|X_{1},\ldots,X_{m}\Big)\Big].\nonumber\\
&\leq&\frac{1}{n}E\Big[E\Big(I(X_{1},Y_{1},F_{m},F)^{2}|X_{1},\ldots,X_{m}\Big)\Big]\nonumber\\
&=&\frac{1}{n}E\Big[E\Big(|I(X_{1},Y_{1},F_{m},F)||X_{1},\ldots,X_{m}\Big)\Big].\nonumber
\end{eqnarray}
Similar to (\ref{a2}), by applying the Chebyshev's inequality, we have
\begin{eqnarray}
W_{nm2}&=&O_{p}\Big(m^{-1/4}n^{-1/2}\Big).\label{a4}
\end{eqnarray}
Therefore, we  obtain  \eqref{a0} and \eqref{a00}.  It is noted that  \eqref{a0} has been proved by Zuo and He \cite{zuo2006} in Lemma 1 (iii). 

\begin{lemma}\label{lem 4}
Under A1 and A3 with $\alpha=2$, we have  
\begin{eqnarray}
\iint I(x,y,F_{m},F)d(F_m-F)(x)dG(y)=O_{p}(m^{-1}).\nonumber
\end{eqnarray}
\end{lemma}
\textbf{Proof of Lemma \ref{lem 4}.} First, we calculate the expectation as follows.
\begin{eqnarray}
&&E\iint I(x,y,F_m,F) d(F_m-F)(x)dG(y)\nonumber\\
&=&E\Big[\frac{1}{m}\sum_{i=1}^{m} E_{Y}I(X_{i},Y,F_{m},F)-E_{X,Y}I(X,Y,F_{m},F)\Big]\nonumber\\
&=&0.\nonumber
\end{eqnarray}
and by A1 and A3 with $\alpha=2$, Lemma \ref{lem 2}, it has
\begin{eqnarray}
&&E\Big[\iint I(x,y,F_{m},F)d(F_m-F)(x)dG(y)\Big]^{2}\nonumber\\
&=&E\Big[\frac{1}{m}\sum_{i=1}^{m} E_{Y}I(X_{i},Y,F_{m},F)-E_{X,Y}I(X,Y,F_{m},F)\Big]^{2}\nonumber\\
&=&E\Big\{E\Big[\Big(\frac{1}{m}\sum_{i=1}^{m} E_{Y}I(X_{i},Y,F_{m},F)-E_{X,Y}I(X,Y,F_{m},F)\Big)^{2}\Big|F_{m}\Big]\Big\}.\nonumber\\
&\leq&\frac{C}{m}E\Big\{E\Big[\Big(E_{Y}I(X_{i},Y,F_{m},F)\Big)^2\Big|F_{m}\Big]\Big\}\nonumber\\
&\leq&\frac{C}{m}E\Big\{E\Big[\Big(E_{Y}|I(X_{i},Y,F_{m},F)|\Big)^2\Big|F_{m}\Big]\Big\}\nonumber\\
&\leq&\frac{C}{m}E\Big\{E\Big[\Big(\sup\limits_{x\in R^{d}}(|D(x;F_{m})-D(x;F)|)\Big)^2\Big|F_{m}\Big]\Big\}=O(m^{-2})
.\nonumber
\end{eqnarray}
Thus, by using the Chebyshev's inequality, it has
\begin{eqnarray}
\iint I(x,y,F_m,F) d(F_m-F)(x)dG(y)=O_{p}(m^{-1}).\nonumber
\end{eqnarray}
Therefore, the proof is complete.

Zuo and He (\cite[Lemma 1(ii)]{zuo2006}) based on conditions A1 and A2, obtain $\iint I(x,y,F_{m},F)d(F_m-F)(x)dG(y)=o_{p}(m^{-1/2})$. In fact, if we choose A1 and A3 with $\alpha=2$ in Lemma \ref{lem 4}, the rate $o_{p}(m^{-1/2})$ can be improved to $O_{p}(m^{-1})$.

\begin{lemma}[Lemma 2 \cite{zuo2006}] \label{lem 5}
Under  A4, we have
$$\iint I(D(x;F_{m})=D(y;F_{m}))dF(x)dG(y)=o(\rho_{m}).$$
\end{lemma}

Lemma \ref{lem 5} relies on A4, and the details can be found in Lemma 2 of Zuo and He \cite{zuo2006}.
\begin{lemma}[\cite{zuo2006}] \label{lem 6}
Under A4, it also has
$$\iint I(x,y,F_{m},F)dF(x)dG(y)=o(\rho_{m}).$$
\end{lemma}

Lemma \ref{lem 6} primarily relies on condition of Lemma \ref{lem 5}. The foundation of the proof can be referenced in Zuo and He \cite{zuo2006}.

\begin{lemma}[Lemma 1\cite{zuo2006}] \label{lem 7}
 There is
$$\iint I(x,y,F)d(F_{m}-F)(x)d(G_{n}-G)(y)=O_{p}\Big(n^{-1/2}m^{-1/2}\Big).$$
\end{lemma}
\textbf{Proof of Lemma \ref{lem 7}.} Let $Z_{nm}=\iint I(x,y,F)d(F_{m}-F)(x)d(G_{n}-G)(y)$
\begin{eqnarray}
EZ_{nm}&=&E\iint I(x,y,F)d(F_{m}-F)(x)d(G_{n}-G)(y)\nonumber\\
&=&\frac{1}{m}\sum_{i=1}^{m}E\Big[\int I(X_{i},y,F)d(G_{n}-G)(y)-E_{X}\Big(\int I(X,y,F)d(G_{n}-G)(y)\Big)\Big]\nonumber\\
&=&0.\nonumber
\end{eqnarray}
Based on the independence of the two samples, the independence within the samples and the boundedness of $I(x, y, F)$, and according to the Cauchy-Schwartz inequality, we obtain that 
\begin{eqnarray}
EZ^{2}_{nm}&=&E\Big[\iint I(x,y,F)d(F_{m}-F)(x)d(G_{n}-G)(y)\Big]^{2}\nonumber\\
&=&E\Big[\frac{1}{n}\sum_{j=1}^{n}\int I(x,Y_{j},F)d(F_{m}-F)(x)-E_{Y}\Big(\int I(x,Y,F)d(F_{m}-F)(x)\Big)\Big]^{2}\nonumber\\
&=&E\Big[\text{Var}\Big(\frac{1}{n}\Big(\sum_{j=1}^{n}\int I(x,Y_{j},F)d(F_{m}-F)(x)\Big)|X_{1},\ldots,X_{m}\Big)\Big].\nonumber\\
&\leq&\frac{1}{n}E\Big[\int I(x,Y_{1},F)d(F_{m}-F)(x)\Big]^{2}\nonumber\\
&=&\frac{1}{n}E\Big[\text{Var}\Big(\frac{1}{m}\sum_{i=1}^{m}I(X_{i},Y_{1},F)|Y_{1}\Big)\Big]\nonumber\\
&\leq&\frac{1}{nm}E\Big[E\Big(\Big(I(X_{1},Y_{1},F)\Big)^{2}|Y_{1}\Big)\Big]\nonumber\\
&=&O(n^{-1}m^{-1}).\nonumber
\end{eqnarray}
Thus, we have
\begin{eqnarray}
Z_{nm}&=&O_{p}\Big(n^{-1/2}m^{-1/2}\Big).\label{l0}
\end{eqnarray}

Lemma \ref{lem 7} was first stated by Zuo and He \cite{zuo2006} in Lemma 1 (i), but the proof was not provided. We provide a detailed proof that provides completeness to the theorem.

\begin{lemma} \label{lem 8}
Under   A6, for $j\neq l$, we have
\begin{eqnarray}
E\Big(I(X_{i},Y_{j},F_{m},F)E_{X}I(X,Y_{l},F_{m},F)\Big)=E\Big(E_{X}I(X,Y_{j},F_{m},F)E_{X}I(X,Y_{l},F_{m},F)\Big)
=0;\label{r4}
\end{eqnarray}
for $i\neq k$,
\begin{eqnarray}
E\Big(I(X_{i},Y_{j},F_{m},F)E_{Y}I(X_{k},Y,F_{m},F)\Big)=E\Big(E_{Y}I(X_{i},Y,F_{m},F)E_{Y}I(X_{k},Y,F_{m},F)\Big)=0;\label{r5}
\end{eqnarray}
for $i\neq k,j\neq l$,
\begin{eqnarray}
E\Big(I(X_{i},Y_{j},F_{m},F)I(X_{k},Y_{l},F_{m},F)\Big)=0.\label{r6}
\end{eqnarray}
for $i=k$ and $j\neq l$, 
\begin{eqnarray}
&&E\Big(I(X_{i},Y_{j},F_{m},F)I(X_{k},Y_{l},F_{m},F)\Big)\nonumber\\
&=&E\Big(I(X_{i},Y_{j},F_{m},F)E_{Y_{l}}I(X_{k},Y_{l},F_{m},F)\Big)\nonumber\\
&=&E\Big(E_{Y_{j}}I(X_{i},Y_{j},F_{m},F)E_{Y_{l}}I(X_{k},Y_{l},F_{m},F)\Big);\label{r7}
\end{eqnarray}
and, for $i\neq k$ and $j=l$,
\begin{eqnarray}
&&E\Big(I(X_{i},Y_{j},F_{m},F)I(X_{k},Y_{l},F_{m},F)\Big)\nonumber\\
&=&E\Big(I(X_{i},Y_{j},F_{m},F)E_{X_{k}}I(X_{k},Y_{l},F_{m},F)\Big)\nonumber\\
&=&E\Big(E_{X_{i}}I(X_{i},Y_{j},F_{m},F)E_{X_{k}}I(X_{k},Y_{l},F_{m},F)\Big).\label{r8}
\end{eqnarray}
\end{lemma}
\textbf{Proof of Lemma \ref{lem 8}.} When 
$j\neq l$, based on the independence of $Y_{j}$ and $Y_{l}$ and A6, we have
\begin{eqnarray}
&&E\Big(E_{X}I(X,Y_{j},F_{m},F)E_{X}I(X,Y_{l},F_{m},F)\Big)\nonumber\\
&=&E\Big\{E\Big[\Big(I(X_{i},Y_{j},F_{m},F)E_{X}I(X,Y_{l},F_{m},F)\Big)\Big|F_{m},Y_j, Y_{l}\Big]\Big\}\nonumber\\
&=&E\Big(I(X_{i},Y_{j},F_{m},F)E_XI(X,Y_{l},F_{m},F)\Big)\nonumber\\
&=&E\Big\{E\Big[\Big(I(X_{i},Y_{j},F_{m},F)E_XI(X,Y_{l},F_{m},F)\Big)\Big|F_{m}\Big]\Big\}\nonumber\\
&=&E\Big\{E_{X_i, Y_j}\Big[\Big(I(X_{i},Y_{j},F_{m},F)\Big)\Big|F_{m}\Big]E_{Y_l}\Big[E_{X}\Big(I(X,Y_{l},F_{m},F)\Big)\Big|F_{m}\Big]\Big\}\nonumber\\
&=&0.\nonumber
\end{eqnarray}
Therefore, 
(\ref{r4}) holds. Similarly, when $i\neq k$ and $i\neq k,j\neq l$, (\ref{r5}) and (\ref{r6}) also hold respectively.
When $i=k$ and $j\neq l$, it has
\begin{eqnarray}
&&E\Big(I(X_{i},Y_{j},F_{m},F)I(X_{k},Y_{l},F_{m},F)\Big)\nonumber\\
&=&E\Big\{E\Big[\Big(I(X_{i},Y_{j},F_{m},F)I(X_{i},Y_{l},F_{m},F)\Big)\Big|X_{i},Y_{j},F_{m}\Big]\Big\}\nonumber\\
&=&E\Big\{I(X_{i},Y_{j},F_{m},F)E_{Y_{l}}\Big[I(X_{i},Y_{l},F_{m},F)\Big|X_{i},Y_{j},F_{m}\Big]\Big\}\nonumber\\
&=&E\Big(I(X_{i},Y_{j},F_{m},F)E_{Y_{l}}I(X_{k},Y_{l},F_{m},F)\Big),\nonumber 
\end{eqnarray}
and
\begin{eqnarray}
&&E\Big(I(X_{i},Y_{j},F_{m},F)I(X_{k},Y_{l},F_{m},F)\Big)\nonumber\\
&=&E\Big\{E\Big[\Big(I(X_{i},Y_{j},F_{m},F)I(X_{i},Y_{l},F_{m},F)\Big)\Big|X_{i},F_{m}\Big]\Big\}\nonumber\\
&=&E\Big\{E_{Y_{j}}\Big[I(X_{i},Y_{j},F_{m},F)\Big|X_{i},F_{m}\Big]E_{Y_{l}}\Big[I(X_{i},Y_{l},F_{m},F)\Big|X_{i},F_{m}\Big]\Big\}\nonumber\\
&=&E\Big(E_{Y_{j}}I(X_{i},Y_{j},F_{m},F)E_{Y_{l}}I(X_{k},Y_{l},F_{m},F)\Big)\nonumber 
\end{eqnarray}
Thus, 
(\ref{r7}) holds. Likewise, when $i\neq k$ and $j=l$, (\ref{r8}) also holds.\\

\begin{lemma}  \label{lem 9}
Under the null hypothesis $F = G$, then by    A3 with $\alpha=2$, A5 and A6, there is
\begin{eqnarray}
&&\iint I(x,y,F_{m},G_{n})d(F_{m}-F)(x)d(G_{n}-G)(y)=O_{p}(m^{-1}n^{-1/2}). \label{f00}
\end{eqnarray}
\end{lemma}
\textbf{Proof of Lemma \ref{lem 9}.}
By expressing the integral as a summation, we have that
\begin{eqnarray}
&&\iint I(x,y,F_{m},F)d(F_{m}-F)(x)d(G_{n}-G)(y)\nonumber\\
&=&\frac{1}{mn}\sum_{i=1}^{m}\sum_{j=1}^{n}I(X_{i},Y_{j},F_{m},F)-\frac{1}{m}\sum_{i=1}^{m}E_{Y}I(X_{i},Y,F_{m},F)\nonumber\\
&&-\frac{1}{n}\sum_{j=1}^{n}E_{X}I(X,Y_{j},F_{m},F)+E_{X,Y}I(X,Y,F_{m},F).\nonumber
\end{eqnarray}
Then, using the law of double expectation, we obtain
\begin{eqnarray}
&&E\iint I(x,y,F_{m},F)d(F_{m}-F)(x)d(G_{n}-G)(y)\nonumber\\
&=&E\Big[\frac{1}{mn}\sum_{i=1}^{m}\sum_{j=1}^{n}I(X_{i},Y_{j},F_{m},F)-\frac{1}{m}\sum_{i=1}^{m}E_{Y}I(X_{i},Y,F_{m},F)\Big]\nonumber\\
&&-E\Big[\frac{1}{n}\sum_{j=1}^{n}E_{X}I(X,Y_{j},F_{m},F)-E_{X,Y}I(X,Y,F_{m},F)\Big]\nonumber\\
&=&E\Big[\frac{1}{nm}\sum_{i=1}^{m}\sum_{j=1}^{n}\Big(I(X_{i},Y_{j},F_{m},F)-E_{Y}I(X_{i},Y,F_{m},F)\Big)\Big]\nonumber\\
&&-E\Big[\frac{1}{n}\sum_{j=1}^{n}\Big(E_{X}I(X,Y_{j},F_{m},F)-E_{X,Y}I(X,Y,F_{m},F)\Big)\Big]\nonumber\\
&=&0,\label{p5}
\end{eqnarray}
and
\begin{eqnarray}
&&\text{Var}\Big(\iint I(x,y,F_{m},F)d(F_{m}-F)(x)d(G_{n}-G)(y)\Big)\nonumber\\
&=&\text{Cov}\Big(\frac{1}{mn}\sum_{i=1}^{m}\sum_{j=1}^{n}I(X_{i},Y_{j},F_{m},F),
\frac{1}{mn}\sum_{k=1}^{m}\sum_{l=1}^{n}I(X_{k},Y_{l},F_{m},F)\Big)\nonumber\\
&&-2\text{Cov}\Big(\frac{1}{mn}\sum_{i=1}^{m}\sum_{j=1}^{n}I(X_{i},Y_{j},F_{m},F),
\frac{1}{m}\sum_{k=1}^{m}E_{Y}I(X_{k},Y,F_{m},F)\Big)\nonumber\\
&&-2\text{Cov}\Big(\frac{1}{mn}\sum_{i=1}^{m}\sum_{j=1}^{n}I(X_{i},Y_{j},F_{m},F),\frac{1}{n}\sum_{l=1}^{n}E_{X}I(X,Y_{l},F_{m},F)
\Big)\nonumber\\
&&+\text{Cov}\Big(\frac{1}{m}\sum_{i=1}^{m}E_{Y}I(X_{i},Y,F_{m},F),
\frac{1}{m}\sum_{k=1}^{m}E_{Y}I(X_{k},Y,F_{m},F)\Big)\nonumber\\
&&+\text{Cov}\Big(\frac{1}{n}\sum_{j=1}^{n}E_{X}I(X,Y_{j},F_{m},F),
\frac{1}{n}\sum_{l=1}^{n}E_{X}I(X,Y_{l},F_{m},F)\Big)\nonumber\\
&&+\text{Cov}\Big(E_{X,Y}I(X,Y,F_{m},F), E_{X,Y}I(X,Y,F_{m},F)\Big)\nonumber\\
&=:&A_{mn1}-2A_{mn2}-2A_{mn3}+A_{mn4}+A_{mn5}+A_{mn6}.\nonumber
\end{eqnarray}
For the first term, we have that
\begin{eqnarray}
A_{mn1}&=&\text{Cov}\Big(\frac{1}{mn}\sum_{i=1}^{m}\sum_{j=1}^{n}I(X_{i},Y_{j},F_{m},F),
\frac{1}{mn}\sum_{k=1}^{m}\sum_{l=1}^{n}I(X_{k},Y_{l},F_{m},F)\Big)\nonumber\\
&=&\frac{1}{m^{2}n^{2}}\sum_{i=k}\sum_{j=l}\text{Cov}\Big(I(X_{i},Y_{j},F_{m},F), I(X_{k},Y_{l},F_{m},F)\Big)\nonumber\\
&&+\frac{1}{m^{2}n^{2}}\sum_{i=k}\sum_{j\neq
l}\text{Cov}\Big(I(X_{i},Y_{j},F_{m},F),I(X_{k},Y_{l},F_{m},F)\Big)\nonumber\\
&&+\frac{1}{m^{2}n^{2}}\sum_{i\neq k}\sum_{j=l}\text{Cov}\Big(I(X_{i},Y_{j},F_{m},F),
I(X_{k},Y_{l},F_{m},F)\Big)\nonumber\\
&&+\frac{1}{m^{2}n^{2}}\sum_{i\neq k}\sum_{j\neq l}\text{Cov}\Big(I(X_{i},Y_{j},F_{m},F),
I(X_{k},Y_{l},F_{m},F)\Big)\nonumber\\
&=&B_{mn1}+B_{mn2}+B_{mn3}+B_{mn4},\nonumber
\end{eqnarray}
where
\begin{eqnarray}
B_{mn1}&=&\frac{1}{m^{2}n^{2}}\sum_{i=k}\sum_{j=l}\text{Cov}\Big(I(X_{i},Y_{j},F_{m},F), I(X_{k},Y_{l},F_{m},F)\Big)\nonumber\\
&=&\frac{1}{mn}\text{Cov}\Big(I(X_{1},Y_{1},F_{m},F), I(X_{1},Y_{1},F_{m},F)\Big)\nonumber\\
&=&\frac{1}{mn}\Big[E\Big(I(X_{1},Y_{1},F_{m},F)\Big)^{2}- \Big(EI(X_{1},Y_{1},F_{m},F)\Big)^{2}\Big].\nonumber
\end{eqnarray}
By  A6, we have
\begin{eqnarray}
E\Big(I(X_{1},Y_{1},F_{m},F)\Big)=0,\label{u1}
\end{eqnarray}
and by (\ref{z1})
\begin{eqnarray}
&&E\Big(I(X_{1},Y_{1},F_{m},F)\Big)^{2}\nonumber\\
&=&E\Big|I(D(X_{1};F_{m})\leq D(Y_{1};F_{m}))-I(D(X_{1};F)\leq D(Y_{1};F))\Big|\nonumber\\
&=&E\Big[P_{X_{1}}(D(X_{1};F_{m})\leq D(Y_{1};F_{m}), D(X_{1};F)> D(Y_{1};F))\nonumber\\
&&+P_{X_{1}}(D(X_{1};F_{m})> D(Y_{1};F_{m}), D(X_{1};F)\leq D(Y_{1};F))\Big]\nonumber\\
&=&E\Big[P_{X_{1}}\Big(D(X_{1};F)\leq D(Y_{1};F)+D(Y_{1};F_{m})-D(Y_{1};F)-D(X_{1};F_{m})+D(X_{1};F),\nonumber\\
&& D(X_{1};F)> D(Y_{1};F)\Big)
\Big]\nonumber\\
&&+E\Big[P_{X_{1}}\Big(D(X_{1};F)>D(Y_{1};F)+D(Y_{1};F_{m})-D(Y_{1};F)-D(X_{1};F_{m})+D(X_{1};F),\nonumber\\
&& D(X_{1};F)\leq D(Y_{1};F)\Big)
\Big]\nonumber\\
&=&T_{1}+T_{2}.\nonumber
\end{eqnarray}
According to (\ref{z1}) and A5, it has
\begin{eqnarray}
T_{1}&=&E\Big[F_{D(X_{1};F)}\Big(D(Y_{1};F)\Big)+f_{D(X_{1};F)}\Big(D(Y_{1};F)\Big)\eta_{1}(Y_{1};F_{m},F)\Big]\nonumber\\
&&-E\Big[F_{D(X_{1};F)}\Big(D(Y_{1};F)\Big)\Big]+O(m^{-1})\nonumber\\
&=&E\Big[f_{D(X_{1};F)}\Big(D(Y_{1};F)\Big)\eta_{1}(Y_{1};F_{m},F)\Big]+O(m^{-1}) \nonumber\\
&=&O(m^{-1}). ~(by~ A5)~~\nonumber
\end{eqnarray}
The proof of $T_{2}$ is similar so it is omitted. Thus,
\begin{eqnarray}
E\Big(I(X_{1},Y_{1},F_{m},G_{n})\Big)^{2}=O(m^{-1}),\nonumber
\end{eqnarray}
So,
\begin{eqnarray}
B_{mn1}=O_{p}(m^{-2}n^{-1}).\label{p9}
\end{eqnarray}
and
\begin{eqnarray}
B_{mn2}&=&\frac{1}{m^{2}n^{2}}\sum_{i=k}\sum_{j\neq
l}\text{Cov}\Big(I(X_{i},Y_{j},F_{m},F),I(X_{k},Y_{l},F_{m},F)\Big)\nonumber\\
&=&\frac{1}{m^{2}n^{2}}\sum_{i=k}\sum_{j\neq
l}\Big[E\Big(I(X_{i},Y_{j},F_{m},F)I(X_{k},Y_{l},F_{m},F)\Big)\nonumber\\
&&-EI(X_{i},Y_{j},F_{m},F)EI(X_{k},Y_{l},F_{m},F)\Big]\nonumber\\
&=&\frac{1}{mn^{2}}\sum_{j\neq
l}\Big[E\Big(E_{Y_{j}}I(X_{1},Y_{j},F_{m},F)E_{Y_{l}}I(X_{1},Y_{l},F_{m},F)\Big)\nonumber~(by~ (\ref{r7})~ in~ Lemma~\ref{lem 8} ~and~(\ref{u1}))\\
&=&\frac{n-1}{mn}E\Big(E_{Y_{1}}I(X_{1},Y_{1},F_{m},F)E_{Y_{2}}I(X_{1},Y_{2},F_{m},F)\Big)
,\label{p10}
\end{eqnarray}
where $B_{mn2}$ will be treated in conjunction with the following $C_{mn1}$ as shown in \ref{p14}.
Similar to (\ref{p10}), we have
\begin{eqnarray}
B_{mn3}&=&\frac{1}{m^{2}n^{2}}\sum_{i\neq k}\sum_{j=l}\text{Cov}\Big(I(X_{i},Y_{j},F_{m},F),
I(X_{k},Y_{l},F_{m},F)\Big)\nonumber\\
&=&\frac{1}{m^{2}n^{2}}\sum_{i\neq
k}\sum_{j=l}\Big[E\Big(I(X_{i},Y_{j},F_{m},F)I(X_{k},Y_{l},F_{m},F)\Big)\nonumber\\
&&-EI(X_{i},Y_{j},F_{m},F)EI(X_{k},Y_{l},F_{m},F)\Big]\nonumber\\
&=&\frac{1}{m^{2}n}\sum_{i\neq
k}\Big[E\Big(E_{X_{i}}I(X_{i},Y_{1},F_{m},F)E_{X_{k}}I(X_{k},Y_{1},F_{m},F)\Big)~(by~ (\ref{r8})~ in~ Lemma~ \ref{lem 8}~and~(\ref{u1}))\nonumber\\
&=&\frac{m-1}{mn}E\Big(E_{X_{1}}I(X_{1},Y_{1},F_{m},F)E_{X_{2}}I(X_{2},Y_{1},F_{m},F)\Big)
,\nonumber
\end{eqnarray}
where $B_{mn3}$ will also be treated in conjunction with $D_{mn1}$ below (see (\ref{p20})). By (\ref{r6}) and (\ref{u1}),
\begin{eqnarray}
B_{mn4}&=&\frac{1}{m^{2}n^{2}}\sum_{i\neq k}\sum_{j\neq l}\text{Cov}\Big(I(X_{i},Y_{j},F_{m},F),
I(X_{k},Y_{l},F_{m},F)\Big)\nonumber\\
&=&\frac{1}{m^{2}n^{2}}\sum_{i\neq k}\sum_{j\neq l}\Big[EI(X_{i},Y_{j},F_{m},F)I(X_{k},Y_{l},F_{m},F)\nonumber\\
&&-EI(X_{i},Y_{j},F_{m},F)EI(X_{k},Y_{l},F_{m},F)\Big]=0
.\nonumber
\end{eqnarray}
Secondly, there is
\begin{eqnarray}
A_{mn2}&=&\text{Cov}\Big(\frac{1}{mn}\sum_{i=1}^{m}\sum_{j=1}^{n}I(X_{i},Y_{j},F_{m},F),
\frac{1}{m}\sum_{k=1}^{m}E_{Y}I(X_{k},Y,F_{m},F)\Big)\nonumber\\
&=&\frac{1}{m^{2}}\sum_{i=k}\text{Cov}\Big(I(X_{i},Y_{1},F_{m},F), E_{Y}I(X_{k},Y,F_{m},F)\Big)\nonumber\\
&&+\frac{1}{m^{2}}\sum_{i\neq k}\text{Cov}\Big(I(X_{i},Y_{1},F_{m},F), E_{Y}I(X_{k},Y,F_{m},F)\Big)\nonumber\\
&=&C_{mn1}+C_{mn2},\nonumber
\end{eqnarray}
where
\begin{eqnarray}
C_{mn1}&=&\frac{1}{m^{2}}\sum_{i=k}\text{Cov}\Big(I(X_{i},Y_{1},F_{m},F), E_{y}I(X_{k},y,F_{m},F)\Big)\nonumber\\
&=&\frac{1}{m^{2}}\sum_{i=k}\Big[E\Big(I(X_{i},Y_{1},F_{m},F)E_{Y}I(X_{k},Y,F_{m},F)\Big)\nonumber\\
&&-EI(X_{i},Y_{1},F_{m},F)E\Big(E_{Y}I(X_{k},Y,F_{m},F)\Big)\Big]\nonumber\\
&=&\frac{1}{m}E\Big(E_{Y_{1}}I(X_{1},Y_{1},F_{m},F)E_{Y}I(X_{1},Y,F_{m},F)\Big)~(by~ (\ref{r7})~ in~ Lemma~ \ref{lem 8}~and~(\ref{u1}))\nonumber\\
&=&\frac{n}{n-1}B_{mn2}.\label{p14}
\end{eqnarray}
Under $H_{0}:F = G$, by (\ref{z1}) and   A3 with $\alpha=2$, it holds
\begin{eqnarray}
&&E\Big(E_{Y_{1}}I(X_{1},Y_{1},F_{m},F)E_{Y}I(X_{1},Y,F_{m},F)\Big)\nonumber\\
&=&E\Big\{\Big[-f_{D(Y_{1};F)}\Big(D(X_{1};F)\Big)\eta_{1}(X_{1};F_{m},F)+O_{p}(m^{-1})\Big]
\times\nonumber\\
&&\Big[-f_{D(Y;F)}\Big(D(X_{1};F)\Big)\eta_{1}(X_{1};F_{m},F)+O_{p}(m^{-1})\Big]\Big\}\nonumber\\
&=&E\Big[-f_{D(Y_{1};F)}\Big(D(X_{1};F)\Big)\eta_{1}(X_{1};F_{m},F)+O_{p}(m^{-1})\Big]^{2}=O(m^{-1}).\label{p15}
\end{eqnarray}
Thus, by (\ref{p14}) and (\ref{p15}), it has
\begin{eqnarray}
B_{mn2}-C_{mn1}=O(m^{-2}n^{-1}).\label{p17}
\end{eqnarray}
By (\ref{r5}) in Lemma \ref{lem 8} and A6,
\begin{eqnarray}
C_{mn2}&=&\frac{1}{m^{2}}\sum_{i\neq k}\text{Cov}\Big(I(X_{i},Y_{1},F_{m},F), E_{y}I(X_{k},y,F_{m},F)\Big)\nonumber\\
&=&\frac{1}{m^{2}}\sum_{i\neq k}\Big[E\Big(I(X_{i},Y_{1},F_{m},F)E_{Y}I(X_{k},Y,F_{m},F)\Big)\nonumber\\
&&-EI(X_{i},Y_{1},F_{m},F)E\Big(E_{Y}I(X_{k},Y,F_{m},F)\Big)\Big]=0.\label{p18}
\end{eqnarray}
For $A_{mn3}$,
\begin{eqnarray}
A_{mn3}&=&\text{Cov}\Big(\frac{1}{mn}\sum_{i=1}^{m}\sum_{j=1}^{n}I(X_{i},Y_{j},F_{m},F),\frac{1}{n}\sum_{l=1}^{n}E_{X}I(X,Y_{l},F_{m},F)
\Big)\nonumber\\
&=&\frac{1}{n^{2}}\sum_{j=l}\text{Cov}\Big(I(X_{1},Y_{j},F_{m},F), E_{X}I(X,Y_{l},F_{m},F)\Big)\nonumber\\
&&+\frac{1}{n^{2}}\sum_{j\neq l}\text{Cov}\Big(I(X_{1},Y_{j},F_{m},F), E_{X}I(X,Y_{l},F_{m},F)\Big)\nonumber\\
&=&D_{mn1}+D_{mn2},\nonumber
\end{eqnarray}
where
\begin{eqnarray}
D_{mn1}&=&\frac{1}{n^{2}}\sum_{j=l}\text{Cov}\Big(I(X_{1},Y_{j},F_{m},F), E_{X}I(X,Y_{l},F_{m},F)\Big)\nonumber\\
&=&\frac{1}{n^{2}}\sum_{j=l}\Big[E\Big(I(X_{1},Y_{j},F_{m},F)E_{X}I(X,Y_{l},F_{m},F)\Big)\nonumber\\
&&-EI(X_{1},Y_{j},F_{m},F)E\Big(E_{X}I(X,Y_{l},F_{m},F)\Big)\Big]\nonumber\\
&=&\frac{1}{n}E\Big(E_{X_{1}}I(X_{1},Y_{1},F_{m},F)E_{X}I(X,Y_{1},F_{m},F)\Big)\nonumber\\
&=&\frac{m}{m-1}B_{mn3}
. ~(by~ (\ref{r8})~ in~ Lemma~ \ref{lem 8})\label{p20}
\end{eqnarray}
Similar to (\ref{p15}), by by condition A3 with $\alpha=2$, it holds
\begin{eqnarray}
&&E\Big(E_{X_{1}}I(X_{1},Y_{1},F_{m},F)E_{X_{2}}I(X_{2},Y_{2},F_{m},F)\Big)\nonumber\\
&=&E\Big\{\Big[-f_{D(X_{1};F)}\Big(D(Y_{1};F)\Big)\eta_{1}(Y_{1};F_{m},F)+O_{p}(m^{-1})\Big]\times\nonumber\\
&&\Big[-f_{D(X_{2};F)}\Big(D(Y_{1};F)\Big)\eta_{1}(Y_{1};F_{m},F)+O_{p}(m^{-1})\Big]\Big\}\nonumber\\
&=&E\Big[-f_{D(X_{1};F)}\Big(D(Y_{1};F)\Big)\eta_{1}(Y_{1};F_{m},F)+O_{p}(m^{-1})\Big]^{2}=O(m^{-1}).\label{p21}
\end{eqnarray}
Thus, by (\ref{p20}) and (\ref{p21}), it has
\begin{eqnarray}
B_{mn3}-D_{mn1}=O(m^{-2}n^{-1}).\label{p23}
\end{eqnarray}
By (\ref{r4}) in Lemma \ref{lem 8} and  A6,
\begin{eqnarray}
D_{mn2}&=&\frac{1}{n^{2}}\sum_{j\neq l}\text{Cov}\Big(I(X_{1},Y_{j},F_{m},F), E_{X}I(X,Y_{l},F_{m},F)\Big)\nonumber\\
&=&\frac{1}{n^{2}}\sum_{j\neq l}\Big[E\Big(I(X_{1},Y_{j},F_{m},F)E_{X}I(X,Y_{l},F_{m},F)\Big)\nonumber\\
&&-EI(X_{1},Y_{j},F_{m},F)E\Big(E_{X}I(X,Y_{l},F_{m},F)\Big)\Big]=0.\label{p24}
\end{eqnarray}
For $A_{mn4}$,
\begin{eqnarray}
A_{mn4}&=&\text{Cov}\Big(\frac{1}{m}\sum_{i=1}^{m}E_{Y}I(X_{i},Y,F_{m},F),
\frac{1}{m}\sum_{k=1}^{m}E_{Y}I(X_{k},Y,F_{m},F)\Big)\nonumber\\
&=&\frac{1}{m^{2}}\sum_{i=k}\text{Cov}\Big(E_{Y}I(X_{i},Y,F_{m},F), E_{Y}I(X_{k},Y,F_{m},F)\Big)\nonumber\\
&&+\frac{1}{m^{2}}\sum_{i\neq k}\text{Cov}\Big(E_{Y}I(X_{i},Y,F_{m},F), E_{Y}I(X_{k},Y,F_{m},F)\Big)\nonumber\\
&=&E_{mn1}+E_{mn2},\nonumber
\end{eqnarray}
where
\begin{eqnarray}
E_{mn1}&=&\frac{1}{m^{2}}\sum_{i=k}\text{Cov}\Big(E_{Y}I(X_{i},Y,F_{m},F), E_{Y}I(X_{k},Y,F_{m},F)\Big)\nonumber\\
&=&\frac{1}{m^{2}}\sum_{i=k}\Big[E\Big(E_{Y}I(X_{i},Y,F_{m},F)E_{Y}I(X_{k},Y,F_{m},F)\Big)\nonumber\\
&&-E\Big(E_{Y}I(X_{i},Y,F_{m},F)\Big)E\Big(E_{Y}I(X_{k},Y,F_{m},F)\Big)\Big]\nonumber\\
&=&\frac{1}{m}\Big[E\Big(E_{Y}I(X_{1},Y,F_{m},F)E_{Y}I(X_{1},Y,F_{m},F)\Big)-EI(X_{1},Y,F_{m},F)EI(X_{1},Y,F_{m},F)\Big]\nonumber\\
&=&C_{mn1}.~(~by~   A6)\label{p26}
\end{eqnarray}
By (\ref{r5}) in Lemma \ref{lem 8} and A6,
\begin{eqnarray}
E_{mn2}&=&\frac{1}{m^{2}}\sum_{i\neq k}\text{Cov}\Big(E_{Y}I(X_{i},Y,F_{m},F), E_{Y}I(X_{k},Y,F_{m},F)\Big)\nonumber\\
&=&\frac{1}{m^{2}}\sum_{i\neq k}\Big[E\Big(E_{Y}I(X_{i},Y,F_{m},F)E_{Y}I(X_{k},Y,F_{m},F)\Big)\nonumber\\
&&-E\Big(E_{Y}I(X_{i},Y,F_{m},F)\Big)E\Big(E_{Y}I(X_{k},Y,F_{m},F)\Big)\Big]=0
.\label{p27}
\end{eqnarray}
For $A_{mn5}$,
\begin{eqnarray}
A_{mn5}&=&\text{Cov}\Big(\frac{1}{n}\sum_{j=1}^{n}E_{X}I(X,Y_{j},F_{m},F),
\frac{1}{n}\sum_{l=1}^{n}E_{X}I(X,Y_{l},F_{m},F)\Big)\nonumber\\
&=&\frac{1}{n^{2}}\sum_{j=l}\text{Cov}\Big(E_{X}I(X,Y_{j},F_{m},F), E_{X}I(X,Y_{l},F_{m},F)\Big)\nonumber\\
&&+\frac{1}{n^{2}}\sum_{j\neq l}\text{Cov}\Big(E_{X}I(X,Y_{j},F_{m},F), E_{X}I(X,Y_{l},F_{m},F)\Big)\nonumber\\
&=&F_{mn1}+F_{mn2}.\nonumber
\end{eqnarray}
Similar to (\ref{p26}), it has
\begin{eqnarray}
F_{mn1}&=&\frac{1}{n^{2}}\sum_{j=l}\text{Cov}\Big(E_{X}I(X,Y_{j},F_{m},F), E_{X}I(X,Y_{l},F_{m},F)\Big)\nonumber\\
&=&\frac{1}{n^{2}}\sum_{j=l}\Big[E\Big(E_{X}I(X,Y_{j},F_{m},F)E_{X}I(X,Y_{l},F_{m},F)\Big)\nonumber\\
&&-E\Big(E_{X}I(X,Y_{j},F_{m},F)\Big)E\Big(E_{X}I(X,Y_{l},F_{m},F)\Big)\Big]\nonumber\\
&=&\frac{1}{n}\Big[E\Big(E_{X}I(X,Y_{1},F_{m},F)E_{X}I(X,Y_{1},F_{m},F)\Big)-EI(X,Y_{1},F_{m},F)EI(X,Y_{1},F_{m},F)\Big]\nonumber\\
&=&D_{mn1}.~(~by~  A6)\label{p28}
\end{eqnarray}
By (\ref{r4}) in Lemma \ref{lem 8} and  A6,
\begin{eqnarray}
F_{mn2}&=&\frac{1}{n^{2}}\sum_{j\neq l}\text{Cov}\Big(E_{X}I(X,Y_{j},F_{m},F), E_{X}I(X,Y_{l},F_{m},F)\Big)\nonumber\\
&=&\frac{1}{n^{2}}\sum_{j\neq l}\Big[E\Big(E_{X}I(X,Y_{j},F_{m},F)E_{X}I(X,Y_{l},F_{m},F)\Big)\nonumber\\
&&-E\Big(E_{X}I(X,Y_{j},F_{m},F)\Big)E\Big(E_{X}I(X,Y_{l},F_{m},F)\Big)\Big]
=0.\label{p31}
\end{eqnarray}
For $A_{mn6}$, we have
\begin{eqnarray}
E_{X,Y}I(X,Y,F_{m},F)=P(D(X;F_{m})\leq D(Y;F_{m}))-P(D(X;F)\leq D(Y;F))=\frac{1}{2}-\frac{1}{2}=0,\nonumber
\end{eqnarray}
which implies
\begin{eqnarray}
A_{mn6}=0.\label{p32}
\end{eqnarray}
Then, combining this with (\ref{p9}), (\ref{p17}), (\ref{p18}), and (\ref{p23})-(\ref{p32}), one has
\begin{eqnarray}
&&\text{Var}\Big(\iint I(X,Y,F_{m},F)d(F_{m}-F)(X)d(G_{n}-G)(Y)\Big)\nonumber\\
&=:&A_{mn1}-2A_{mn2}-2A_{mn3}+A_{mn4}+A_{mn5}+A_{mn6}\nonumber\\
&=&B_{mn1}+B_{mn2}+B_{mn3}+B_{mn4}-2C_{mn1}-2C_{mn2}-2D_{mn1}-2D_{mn2}\nonumber\\
&&+E_{mn1}+E_{mn2}+F_{mn1}+F_{mn2}+A_{mn6}\nonumber\\
&=&B_{mn1}+(B_{mn2}-C_{mn1})+(B_{mn3}-D_{mn1})+A_{mn6}\nonumber\\
&=&O(m^{-2}n^{-1}).\label{p34}
\end{eqnarray}
Based on (\ref{p5}) and (\ref{p34}), the desired result can be obtained by Chebyshev's inequality.\\

Lemma \ref{lem 9} plays an important role in our proof of the remainder term in Theorem \ref{The 2}, but it has not been mentioned elsewhere in the literature. We obtained this result based mainly on the Hoeffding decomposition and the Cox-Reid extension. The theoretical results obtained will pave the way for the study of $Q$-statistics.   

\begin{lemma}\label{lem 10}
Under   A3 and A5,  we have 
\begin{eqnarray}
\frac{1}{n}\sum\limits_{j=1}^{n}f_{D(X;F)}\Big(D(Y_{j};F)\Big)\eta_{1}(Y_{j};F_{m},F)=O_{p}(n^{-1/2}m^{-1/2}),\label{t0}
\end{eqnarray}
and
\begin{eqnarray}
\frac{1}{m}\sum\limits_{i=1}^{m}f_{D(Y;F)}\Big(D(X_{i};F)\Big)\eta_{1}(X_{i};F_{m},F)=O_{p}(m^{-1}).\label{t1}
\end{eqnarray}
\end{lemma}
\textbf{Proof of Lemma \ref{lem 10}.} By A5, we obtain
\begin{eqnarray}
E\Big[\frac{1}{n}\sum\limits_{j=1}^{n}f_{D(X;F)}\Big(D(Y_{j};F)\Big)\eta_{1}(Y_{j};F_{m},F)\Big]=0.\nonumber
\end{eqnarray}
Given the independence and identical distribution among $Y_{1},\cdots, Y_{n}$ and their independence from $F_{m}$, we can conclude that
\begin{eqnarray}
&&E\Big[\frac{1}{n}\sum\limits_{j=1}^{n}f_{D(X;F)}\Big(D(Y_{j};F)\Big)\eta_{1}(Y_{j};F_{m},F)\Big]^2\nonumber\\
&=&\frac{1}{n}E\Big[f_{D(X;F)}\Big(D(Y_{1};F)\Big)\eta_{1}(Y_{1};F_{m},F)\Big]^2\nonumber\\
&=&O(n^{-1}m^{-1}).~(by~   A3~with ~\alpha=2)\nonumber
\end{eqnarray}
The proof of (\ref{t0}) is completed by Chebyshev's inequality. Next, we need to demonstrate that (\ref{t1}) holds.
Similarly, by A5, we have
\begin{eqnarray}
E\Big[\frac{1}{m}\sum\limits_{i=1}^{m}f_{D(Y;F)}\Big(D(X_{i};F)\Big)\eta_{1}(X_{i};F_{m},F)\Big]=0,
\nonumber
\end{eqnarray}
and
\begin{eqnarray}
&&E\Big[\frac{1}{m}\sum\limits_{i=1}^{m}f_{D(Y;F)}\Big(D(X_{i};F)\Big)\eta_{1}(X_{i};F_{m},F)\Big]^2\nonumber\\
&=&\frac{1}{m^2}E\Big[\sum\limits_{i=1}^{m}f_{D(Y;F)}\Big(D(X_{i};F)\Big)\eta_{1}(X_{i};F_{m},F)\Big]^2\nonumber\\
&=&\frac{1}{m^2}\sum\limits_{i=k}E\Big[f_{D(Y;F)}\Big(D(X_{i};F)\Big)f_{D(Y;F)}\Big(D(X_{k};F)\Big)\eta_{1}(X_{i};F_{m},F)\eta_{1}(X_{k};F_{m},F)\Big]\nonumber\\
&&+\frac{1}{m^2}\sum\limits_{i\neq k}E\Big[f_{D(Y;F)}\Big(D(X_{i};F)\Big)f_{D(Y;F)}\Big(D(X_{k};F)\Big)\eta_{1}(X_{i};F_{m},F)\eta_{1}(X_{k};F_{m},F)\Big]\nonumber\\
&=&T_{m1}+T_{m2}.\label{t2}
\end{eqnarray}
For $T_{m1}$, by   A3 with $\alpha=2$, we have
\begin{eqnarray}
T_{m1}=\frac{1}{m}E\Big[f_{D(Y;F)}\Big(D(X_{1};F)\Big)\eta_{1}(X_{1};F_{m},F)\Big]^2=O_{p}(m^{-2}).\label{t3}
\end{eqnarray}
For $T_{m2}$, since $X_{1},\cdots, X_{n}$ are independently and identically distributed. From A5 it follows that
\begin{eqnarray}
T_{m2}&=&\frac{m(m-1)}{m^2}E\Big\{E\Big[f_{D(Y;F)}\Big(D(X_{1};F)\Big)\eta_{1}(X_{1};F_{m},F)f_{D(Y;F)}\Big(D(X_{2};F)\Big)\eta_{1}(X_{2};F_{m},F)\Big|\Gamma_{m}\Big]\Big\}\nonumber\\
&=&\frac{m(m-1)}{m^2}E\Big\{E\Big[f_{D(Y;F)}\Big(D(X_{1};F)\Big)\eta_{1}(X_{1};F_{m},F)\Big|\Gamma_{m}\Big]E\Big[f_{D(Y;F)}\Big(D(X_{2};F)\Big)\eta_{1}(X_{2};F_{m},F)\Big|\Gamma_{m}\Big]\Big\}\nonumber\\
&=&0.\nonumber
\end{eqnarray}
Combining this with (\ref{t2}) and (\ref{t3}), and then using Chebyshev's inequality, we can obtain (\ref{t1}).

Lemma \ref{lem 10} plays an important role in our proof of the main term in Theorem 2.
\subsection{Proof of Theorem}
\textbf{Proof of Theorem \ref{The 1}.} It is obvious to obtain that
\begin{eqnarray}
&&Q(F_{m},G_{n})-Q(F,G)\nonumber\\
&=&\iint I(x,y,F_{m})dF_{m}(x)dG_{n}(y)-\iint I(x,y,F)dF(x)dG(y) \nonumber\\
&=&\iint I(x,y,F)dF(x)d(G_{n}-G)(y) \nonumber\\
&&+\iint I(x,y,F)d(F_{m}-F)(x)dG(y)+\text{ZH}_{mn},\label{f1}
\end{eqnarray}
where 
\begin{eqnarray}
\text{ZH}_{mn}&=&\iint I(x,y,F_{m})dF_{m}(x)dG_{n}(y)-\iint I(x,y,F)dF(x)dG(y)\nonumber\\
&&-\iint I(x,y,F)dF(x)d(G_{n}-G)(y)-\iint I(x,y,F)d(F_{m}-F)(x)dG(y)\nonumber\\
&=&\iint I(x,y,F_{m},F)dF_{m}(x)dG_{n}(y)+\iint I(x,y,F)dF_{m}(x)dG_{n}(y)\nonumber\\
&&-\iint I(x,y,F)dF(x)dG_{n}(y)-\iint I(x,y,F)d(F_{m}-F)(x)dG(y)\nonumber\\
&=&\iint I(x,y,F_{m},F)dF_{m}(x)dG_{n}(y)+\iint I(x,y,F)d(F_{m}-F)(x)dG_{n}(y)\nonumber\\
&&-\iint I(x,y,F)d(F_{m}-F)(x)dG(y)\nonumber\\
&=&\iint I(x,y,F_{m},F)dF_{m}(x)d(G_{n}-G)(y)+\iint I(x,y,F_{m},F)dF_{m}(x)dG(y)\nonumber\\
&&+\iint I(x,y,F)d(F_{m}-F)(x)d(G_{n}-G)(y)\nonumber\\
&=&\iint I(x,y,F_{m},F)dF_{m}(x)d(G_{n}-G)(y)+\iint I(x,y,F_{m},F)d(F_{m}-F)(x)dG(y)\nonumber\\
&&+\iint I(x,y,F_{m},F)dF(x)dG(y)+\iint I(x,y,F)d(F_{m}-F)(x)d(G_{n}-G)(y).\label{f2}
\end{eqnarray}
Based on (\ref{a00}) in Lemma \ref{lem 3}, Lemma \ref{lem 4}, Lemma \ref{lem 6}, Lemma \ref{lem 7} 
and given that $m/n$ tends to a positive constant, we have
\begin{eqnarray}
\text{ZH}_{mn}&=&O_{p}\Big(m^{-1/4}n^{-1/2}\Big)+O_{p}(m^{-1})+o(\rho_{m})+O_{p}\Big(n^{-1/2}m^{-1/2}\Big)\nonumber\\
&=&O_{p}\Big(m^{-3/4}\Big)+o(\rho_{m}).\label{f3}
\end{eqnarray}
Therefore, we complete the proof of Theorem \ref{The 1}.\\
\textbf{Proof of Theorem \ref{The 2}.} By Hoeffding decomposition, it has
\begin{eqnarray}
&&Q(F_{m},G_{n})-Q(F,G)\nonumber\\
&=&\iint I(x,y,F_{m})dF_{m}(x)dG_{n}(y)-\iint I(x,y,F)dF(x)dG(y) \nonumber\\
&=&\Big[\iint I(x,y,F_{m})dF(x)dG_{n}(y)-\iint I(x,y,F)dF(x)dG(y)\Big]\nonumber\\
&&+\Big[\iint I(x,y,F_{m})dF_{m}(x)dG(y)-\iint I(x,y,F)dF(x)dG(y)\Big]+R_{mn}\nonumber\\
&=&L_{mn1}+L_{mn2}+R_{mn},\label{p1}
\end{eqnarray}
where
\begin{eqnarray}
R_{mn}&=&\iint I(x,y,F_{m})dF_{m}(x)dG_{n}(y)-\iint I(x,y,F)dF(x)dG(y)\nonumber\\
&&-\iint I(x,y,F_{m})dF(x)dG_{n}(y)+\iint I(x,y,F)dF(x)dG(y)\nonumber\\
&&-\iint I(x,y,F_{m})dF_{m}(x)dG(y)+\iint I(x,y,F)dF(x)dG(y)\nonumber\\
&=&\iint I(x,y,F_{m})dF_{m}(x)dG_{n}(y)-\iint I(x,y,F_{m})dF(x)dG_{n}(y)\nonumber\\
&&-\iint I(x,y,F_{m})dF_{m}(x)dG(y)+\iint I(x,y,F)dF(x)dG(y).\label{p2}
\end{eqnarray}
For $L_{mn1}$, by (\ref{z1}) and Lemma \ref{lem 10}, we have
\begin{eqnarray}
L_{mn1}&=&\iint I(x,y,F_{m})dF(x)dG_{n}(y)-\iint I(x,y,F)dF(x)dG(y)\nonumber\\
&=&\frac{1}{n}\sum\limits_{j=1}^{n}E_{X}I(X,Y_{j},F_{m})-\frac{1}{2}\nonumber\\
&=&\frac{1}{n}\sum\limits_{j=1}^{n}\Big[P_{X}\Big(D(X;F_{m})\leq D(Y_{j};F_{m})\Big)-\frac{1}{2}\nonumber\\
&=&\frac{1}{n}\sum\limits_{j=1}^{n}P_{X}\Big(D(X;F)\leq D(Y_{j};F)+D(Y_{j};F_{m})-D(Y_{j};F)-D(X;F_{m})+D(X;F)\Big)-\frac{1}{2}\nonumber\\
&=&\frac{1}{n}\sum\limits_{j=1}^{n}\Big[F_{D(X;F)}\Big(D(Y_{j};F)\Big)+f_{D(X;F)}\Big(D(Y_{j};F)\Big)\eta_{1}(Y_{j};F_{m},F)+O_{p}(m^{-1})\Big]-\frac{1}{2}\nonumber\\
&=&\frac{1}{n}\sum\limits_{j=1}^{n}\Big[f_{D(X;F)}\Big(D(Y_{j};F)\Big)\eta_{1}(Y_{j};F_{m},F)+O_{p}(m^{-1})\Big]\nonumber\\
&&+\frac{1}{n}\sum\limits_{j=1}^{n}\Big[F_{D(X;F)}\Big(D(Y_{j};F)\Big)-\frac{1}{2}\Big]
\nonumber\\
&=&\frac{1}{n}\sum\limits_{j=1}^{n}\Big[F_{D(X;F)}\Big(D(Y_{j};F)\Big)-\frac{1}{2}\Big]+O_{p}(m^{-1}).
\label{p3}
\end{eqnarray}
Similarly, for $L_{mn2}$, it has
\begin{eqnarray}
L_{mn2}&=&\iint I(x,y,F_{m})dF_{m}(x)dG(y)-\iint I(x,y,F)dF(x)dG(y)\nonumber\\
&=&\frac{1}{m}\sum\limits_{i=1}^{m}E_{Y}I(X_{i},Y,F_{m})-\frac{1}{2}\nonumber\\
&=&\frac{1}{m}\sum\limits_{i=1}^{m}P_{Y}\Big(D(X_{i};F_{m})\leq D(Y;F_{m})\Big)-\frac{1}{2}\nonumber\\
&=&\frac{1}{m}\sum\limits_{i=1}^{m}P_{Y}\Big(D(X_{i};F)\leq D(Y;F)+D(Y;F_{m})-D(Y;F)-D(X_{i};F_{m})+D(X_{i};F)\Big)-\frac{1}{2}\nonumber\\
&=&\frac{1}{m}\sum\limits_{i=1}^{m}\Big[1-F_{D(Y;F)}(D(X_{i};F))-f_{D(Y;F)}(D(X_{i};F))\eta_{1}(X_{i};F_{m},F)+O_{p}(m^{-1})\Big]
-\frac{1}{2}\nonumber\\
&=&\frac{1}{m}\sum\limits_{i=1}^{m}\Big[-f_{D(Y;F)}(D(X_{i};F))\eta_{1}(X_{i};F_{m},F)+O_{p}(m^{-1})\Big]\nonumber\\
&&+\frac{1}{m}\sum\limits_{i=1}^{m}\Big[\frac{1}{2}-F_{D(Y;F)}\Big(D(X_{i};F)\Big)\Big]
\nonumber\\
&=&\frac{1}{m}\sum\limits_{i=1}^{m}\Big[\frac{1}{2}-F_{D(Y;F)}\Big(D(X_{i};F)\Big)\Big]+O_{p}(m^{-1}).
\label{pp3}
\end{eqnarray}
For $R_{mn}$, there is
\begin{eqnarray}
R_{mn}
&=&\iint I(x,y,F_{m})dF_{m}(x)dG_{n}(y)-\iint I(x,y,F_{m})dF(x)dG_{n}(y)\nonumber\\
&&-\iint I(x,y,F_{m})dF_{m}(x)dG(y)+\iint I(x,y,F)dF(x)dG(y).\nonumber\\
&=&\iint I(x,y,F_{m},F)d(F_{m}-F)(x)d(G_{n}-G)(y)-\iint I(x,y,F_{m})dF(x)dG(y)\nonumber\\
&&+\iint I(x,y,F)d(F_{m}-F)(x)d(G_{n}-G)(y)+\iint I(x,y,F)dF(x)dG(y)\nonumber\\
&=&-\iint I(x,y,F_{m},F)dF(x)dG(y)\nonumber\\
&&+\iint I(x,y,F)d(F_{m}-F)(x)d(G_{n}-G)(y)\nonumber\\
&&+\iint I(x,y,F_{m},F)d(F_{m}-F)(x)d(G_{n}-G)(y).\nonumber
\end{eqnarray}
 Using Lemma \ref{lem 6}, Lemma \ref{lem 7}, and Lemma \ref{lem 9}, we have
\begin{eqnarray}
R_{mn}&=&o(\rho_{m})+O_{p}\Big(n^{-1/2}m^{-1/2}\Big)+O_{p}\Big(m^{-1}n^{-1/2}\Big)\nonumber\\
&=&o(\rho_{m})+O_{p}(m^{-1}).\label{pp5}
\end{eqnarray}
Therefore, we obtain
\begin{eqnarray}
&&Q(F_{m},G_{n})-Q(F,G)\nonumber\\
&=&\frac{1}{n}\sum\limits_{j=1}^{n}\Big[F_{D(X;F)}\Big(D(Y_{j};F)\Big)-\frac{1}{2}\Big]+\frac{1}{m}\sum\limits_{i=1}^{m}\Big[\frac{1}{2}-F_{D(Y;F)}\Big(D(X_{i};F)\Big)\Big]\nonumber\\
&&+o(\rho_{m})+O_{p}(m^{-1})
.\label{pppp5}
\end{eqnarray}
The proof of (\ref{ww1}) is completed by (\ref{p1}), (\ref{p3})-(\ref{pp5}). 
Similarly, 
\begin{eqnarray}
&&Q(G_{n},F_{m})-Q(G,F)\nonumber\\
&=&\iint I(y,x,G_{n})dG_{n}(y)dF_{m}(x)-\iint I(y,x,F)dF(x)dG(y) \nonumber\\
&=&\Big[\iint I(y,x,G_{n})dG_{n}(y)dF(x)-\iint I(y,x,F)dF(x)dG(y)\Big]\nonumber\\
&&+\Big[\iint I(y,x,G_{n})dF_{m}(x)dG(y)-\iint I(y,x,F)dF(x)dG(y)\Big]+\tilde{R}_{mn} \nonumber\\
&=&L^{'}_{mn1}+L^{'}_{mn}+\tilde{R}_{mn},\nonumber
\end{eqnarray}
where
\begin{eqnarray}
\tilde{R}_{mn}&=&\iint I(y,x,G_{n})dF_{m}(x)dG_{n}(y)-\iint I(y,x,F)dF(x)dG(y)\nonumber\\
&&-\iint I(y,x,G_{n})dF(x)dG_{n}(y)+\iint I(y,x,F)dF(x)dG(y)\nonumber\\
&&-\iint I(y,x,G_{n})dF_{m}(x)dG(y)+\iint I(y,x,F)dF(x)dG(y)\nonumber\\
&=&\iint I(y,x,G_{n})dF_{m}(x)dG_{n}(y)-\iint I(y,x,G_{n})dF(x)dG_{n}(y)\nonumber\\
&&-\iint I(y,x,G_{n})dF_{m}(x)dG(y)+\iint I(y,x,F)dF(x)dG(y).\nonumber
\end{eqnarray}
By definition of $I(x,y,H)$, it has $I(x,y,H)=1-I(y,x,H)$. Similar to the proof of $Q(F_{m},G_{n})$ and $m/n$ tends to a positive constant, we obtain 
\begin{eqnarray}
L^{'}_{mn1}=\frac{1}{n}\sum\limits_{j=1}^{n}\Big[\frac{1}{2}-F_{D(X;F)}\Big(D(Y_{j};F)\Big)\Big]+O_{p}(m^{-1}),
\nonumber
\end{eqnarray}
\begin{eqnarray}
L^{'}_{mn2}=\frac{1}{m}\sum\limits_{i=1}^{m}\Big[F_{D(Y;F)}\Big(D(X_{i};F)\Big)-\frac{1}{2}\Big]+O_{p}(m^{-1}).
\nonumber
\end{eqnarray}
and
\begin{eqnarray}
\tilde{R}_{mn}=O_{p}(m^{-1})+o(\rho_{m}).\nonumber
\end{eqnarray}
Therefore, we obtain 
\begin{eqnarray}
&&Q(G_{n},F_{m})-Q(G,F)\nonumber\\
&=&\frac{1}{n}\sum\limits_{j=1}^{n}\Big[\frac{1}{2}-F_{D(X;F)}\Big(D(Y_{j};F)\Big)\Big]+\frac{1}{m}\sum\limits_{i=1}^{m}\Big[F_{D(Y;F)}\Big(D(X_{i};F)\Big)-\frac{1}{2}\Big]\nonumber\\
&&+o(\rho_{m})+O_{p}(m^{-1}).\nonumber
\end{eqnarray}
Thus, we complete the proof of (\ref{w1}) in Theorem \ref{The 2}.\\
\textbf{Proof of corollary \ref{cor 2}.}
Under the null hypothesis, $F$ and $G$ have the same distribution and $Q(F,G)=Q(G,F)=\frac{1}{2}$. Since $m/n$ tends to a positive constant, by Theorem \ref{The 2} and the Hoeffding decomposition,  there is 
\begin{eqnarray}
M_{m,n}&=&\Big[\frac{1}{12}\Big(\frac{1}{m}+\frac{1}{n}\Big)\Big]^{-1}\max\Big[\Big(Q(F_{m},G_{n})-1/2\Big)^2, \Big(Q(G_{n},F_{m})-1/2\Big)^2\Big] \nonumber\\
&=&\Big[\frac{1}{12}\Big(\frac{1}{m}+\frac{1}{n}\Big)\Big]^{-1}\max\Big\{\frac{1}{n}\sum\limits_{j=1}^{n}\Big[F_{D(X;F)}\Big(D(Y_{j};F)\Big)-\frac{1}{2}\Big]+\frac{1}{m}\sum\limits_{i=1}^{m}\Big[\frac{1}{2}-F_{D(Y;F)}\Big(D(X_{i};F)\Big)\Big]
\nonumber\\
&&+o(\rho_{m})+O_{p}(m^{-1}), \frac{1}{n}\sum\limits_{j=1}^{n}\Big[\frac{1}{2}-F_{D(X;F)}\Big(D(Y_{j};F)\Big)\Big]+\frac{1}{m}\sum\limits_{i=1}^{m}\Big[F_{D(Y;F)}\Big(D(X_{i};F)\Big)-\frac{1}{2}\Big]\nonumber\\
&&+o(\rho_{m})+O_{p}(m^{-1})\Big\}^2
\nonumber\\
&=&\Big[\frac{1}{12}\Big(\frac{1}{m}+\frac{1}{n}\Big)\Big]^{-1}\Big\{\frac{1}{n}\sum\limits_{j=1}^{n}\Big[F_{D(X;F)}\Big(D(Y_{j};F)\Big)-\frac{1}{2}\Big]+\frac{1}{m}\sum\limits_{i=1}^{m}\Big[\frac{1}{2}-F_{D(Y;F)}\Big(D(X_{i};F)\Big)\Big]\Big\}^2\nonumber\\
&&+o_{p}\Big(m^{1/2}\rho_{m}\Big)+ O_{p}(m^{-1/2}),\nonumber
\end{eqnarray}
and
\begin{eqnarray}
M^{*}_{m,n}&=&\Big[\frac{1}{12}\Big(\frac{1}{m}+\frac{1}{n}\Big)\Big]^{-1/2}\Big[\frac{1}{2}-\min\Big(Q(F_{m},G_{n}),Q(G_{n},F_{m})\Big)\Big]\nonumber\\
&=&\Big[\frac{1}{12}\Big(\frac{1}{m}+\frac{1}{n}\Big)\Big]^{-1/2}\Big[-\min\Big(Q(F_{m},G_{n})-\frac{1}{2},Q(G_{n},F_{m})-\frac{1}{2}\Big)\Big]\nonumber\\
&=&\Big[\frac{1}{12}\Big(\frac{1}{m}+\frac{1}{n}\Big)\Big]^{-1/2}\Big|Q(F_{m},G_{n})-\frac{1}{2}\Big|\nonumber\\
&=&\Big[\frac{1}{12}\Big(\frac{1}{m}+\frac{1}{n}\Big)\Big]^{-1/2}\Big|\frac{1}{n}\sum\limits_{j=1}^{n}\Big[F_{D(X;F)}\Big(D(Y_{j};F)\Big)-\frac{1}{2}\Big]+\frac{1}{m}\sum\limits_{i=1}^{m}\Big[\frac{1}{2}-F_{D(Y;F)}\Big(D(X_{i};F)\Big)\Big]
\nonumber\\
&&+o(\rho_{m})+O_{p}(m^{-1})\Big|\nonumber\\
&=&\Big[\frac{1}{12}\Big(\frac{1}{m}+\frac{1}{n}\Big)\Big]^{-1/2}\Big|\frac{1}{n}\sum\limits_{j=1}^{n}\Big[F_{D(X;F)}\Big(D(Y_{j};F)\Big)-\frac{1}{2}\Big]+\frac{1}{m}\sum\limits_{i=1}^{m}\Big[\frac{1}{2}-F_{D(Y;F)}\Big(D(X_{i};F)\Big)\Big]\Big|
\nonumber\\
&&+o\Big(m^{1/2}\rho_{m}\Big)+O_{p}(m^{-1/2}).\nonumber
\end{eqnarray}
\textbf{Proof of corollary \ref{cor 3}.}
With Corollary \ref{cor 2} and the definition of the chi-square distribution, Corollary \ref{cor 3} immediately follows.
\end{document}